\documentclass[11pt]{article}

\pdfoutput=1

\usepackage{textcomp}
\usepackage[pdftex]{graphicx}
\usepackage[cmex10]{amsmath}

\usepackage{algorithmic,algorithm}
\usepackage{array}
\usepackage[caption=false,font=footnotesize]{subfig}
\usepackage{url}

\usepackage[backend=bibtex, doi=false, isbn=false, url=false]{biblatex} 
\makeatletter
\def\blx@bblfile@bibtex{%
  \blx@secinit
  \begingroup
  \blx@bblstart
\begingroup
\makeatletter
\@ifundefined{ver@biblatex.sty}
  {\@latex@error
     {Missing 'biblatex' package}
     {The bibliography requires the 'biblatex' package.}
      \aftergroup }
  {}
\endgroup

\entry{michau_estimating_2015-1}{inproceedings}{}
  \name{author}{6}{}{%
    {{}%
     {Michau}{M.}%
     {Gabriel}{G.}%
     {}{}%
     {}{}}%
    {{}%
     {Pustelnik}{P.}%
     {Nelly}{N.}%
     {}{}%
     {}{}}%
    {{}%
     {Borgnat}{B.}%
     {Pierre}{P.}%
     {}{}%
     {}{}}%
    {{}%
     {Abry}{A.}%
     {Patrice}{P.}%
     {}{}%
     {}{}}%
    {{}%
     {Nantes}{N.}%
     {Alfredo}{A.}%
     {}{}%
     {}{}}%
    {{}%
     {Chung}{C.}%
     {Edward}{E.}%
     {}{}%
     {}{}}%
  }
  \keyw{Bluetooth,car trajectory,Convex optimization,convex optimization
  resolution algorithm,Estimation,inference on network,inverse problem,inverse
  problems,link dependent origin-destination matrix estimation,mobile
  computing,mobile radio,OD matrices,Roads,sample trajectories,Traffic
  counts,traffic engineering computing,Trajectory,trajectory
  information,transport network,urban transport,Vehicles}
  \strng{namehash}{MGPNBP+1}
  \strng{fullhash}{MGPNBPAPNACE1}
  \field{abstract}{%
  In transport networks, Origin-Destination matrices (ODM) are classically
  estimated from road traffic counts whereas recent technologies grant also
  access to sample car trajectories. One example is the deployment in cities of
  Bluetooth scanners that measure the trajectories of Bluetooth equipped cars.
  Exploiting such sample trajectory information, the classical ODM estimation
  problem is here extended into a link-dependent ODM (LODM) one. This much
  larger size estimation problem is formulated here in a variational form as an
  inverse problem. We develop a convex optimization resolution algorithm that
  incorporates network constraints. We study the result of the proposed
  algorithm on simulated network traffic.%
  }
  \field{booktitle}{2015 {{IEEE International Conference}} on {{Acoustics}},
  {{Speech}} and {{Signal Processing}} ({{ICASSP}})}
  \verb{doi}
  \verb 10.1109/ICASSP.2015.7179019
  \endverb
  \field{eventtitle}{2015 IEEE International Conference on Acoustics, Speech
  and Signal Processing (ICASSP)}
  \field{pages}{5480\bibrangedash 5484}
  \field{title}{Estimating Link-Dependent Origin-Destination Matrices from
  Sample Trajectories and Traffic Counts}
  \list{location}{1}{%
    {{Brisbane, Australia}}%
  }
  \verb{file}
  \verb IEEE Xplore Abstract Record:D\:\\QUT\\Bibliography\\Zotero\\storage\\N6
  \verb ZVZIP4\\articleDetails.html:
  \endverb
  \field{day}{19}
  \field{month}{04}
  \field{year}{2015}
  \field{endday}{24}
  \field{endmonth}{04}
  \field{endyear}{2015}
\endentry

\entry{michau_estimating_2015}{inproceedings}{}
  \name{author}{6}{}{%
    {{}%
     {Michau}{M.}%
     {Gabriel}{G.}%
     {}{}%
     {}{}}%
    {{}%
     {Borgnat}{B.}%
     {Pierre}{P.}%
     {}{}%
     {}{}}%
    {{}%
     {Pustelnik}{P.}%
     {Nelly}{N.}%
     {}{}%
     {}{}}%
    {{}%
     {Abry}{A.}%
     {Patrice}{P.}%
     {}{}%
     {}{}}%
    {{}%
     {Nantes}{N.}%
     {Alfredo}{A.}%
     {}{}%
     {}{}}%
    {{}%
     {Chung}{C.}%
     {Edward}{E.}%
     {}{}%
     {}{}}%
  }
  \strng{namehash}{MGBPPN+1}
  \strng{fullhash}{MGBPPNAPNACE1}
  \field{abstract}{%
  Origin-Destination matrices (ODM) estimation can benefits of the availability
  of sample trajectories which can be measured thanks to recent technologies.
  This paper focus on the case of transport networks where traffic counts are
  measured by magnetic loops and sample trajectories available. An example of
  such network is the city of Brisbane, where Bluetooth detectors are now
  operating. This additional data source is used to extend the classical ODM
  estimation to a link-specific ODM (LODM) one using a convex optimisation
  resolution that incorporates networks constraints as well. The proposed
  algorithm is assessed on a simulated network.%
  }
  \field{booktitle}{{{XXV GRETSI}}}
  \field{eventtitle}{XXV GRETSI}
  \field{title}{Estimating link-dependent origin-destination matrices from
  sample trajectories and traffic counts}
  \verb{url}
  \verb http://eprints.qut.edu.au/86449/
  \endverb
  \field{langid}{french}
  \list{location}{1}{%
    {{Lyon, France}}%
  }
  \verb{file}
  \verb Michau et al_2015_Estimating Link-Dependent Origin-Destination Matrices
  \verb  from Sample Trajectories.pdf:D\:\\QUT\\Bibliography\\Zotero\\storage\\
  \verb STQFN3T7\\Michau et al_2015_Estimating Link-Dependent Origin-Destinatio
  \verb n Matrices from Sample Trajectories.pdf:application/pdf
  \endverb
  \field{day}{08}
  \field{month}{09}
  \field{year}{2015}
\endentry

\entry{willumsen_estimation_1978}{article}{}
  \name{author}{1}{}{%
    {{}%
     {Willumsen}{W.}%
     {Luis~G}{L.~G.}%
     {}{}%
     {}{}}%
  }
  \strng{namehash}{WLG1}
  \strng{fullhash}{WLG1}
  \field{shorttitle}{Estimation of an {{OD Matrix}} from {{Traffic Counts–A
  Review}}}
  \field{title}{Estimation of an {{OD}} Matrix from Traffic Counts: A Review}
  \verb{file}
  \verb Willumsen_1978_Estimation of an OD Matrix from Traffic Counts–A Revie
  \verb w.pdf:D\:\\QUT\\Bibliography\\Zotero\\storage\\P2SKGTUC\\Willumsen_1978
  \verb _Estimation of an OD Matrix from Traffic Counts–A Review.pdf:applicat
  \verb ion/pdf
  \endverb
  \field{journaltitle}{Institute of Transport Studies, Universities of Leed}
  \field{year}{1978}
\endentry

\entry{coates_internet_2002}{article}{}
  \name{author}{4}{}{%
    {{}%
     {Coates}{C.}%
     {M.}{M.}%
     {}{}%
     {}{}}%
    {{}%
     {Hero}{H.}%
     {A.O.}{A.}%
     {}{}%
     {}{}}%
    {{}%
     {Nowak}{N.}%
     {R.}{R.}%
     {}{}%
     {}{}}%
    {{}%
     {Yu}{Y.}%
     {Bin}{B.}%
     {}{}%
     {}{}}%
  }
  \keyw{anomalous/malicious behavior detection,array processing,Array signal
  processing,distributed network,dynamic routing,e-commerce,inference
  mechanisms,inference problems,Internet,Internet tomography,inverse
  problems,IP networks,Monitoring,multicast network tomography,network
  monitoring,network servers,network traffic measurements,optimized service
  provision,performance
  evaluation,reviews,routers,Routing,servers,service-level verification,signal
  processing,Signal processing algorithms,signal processing problems,signal
  processing theory,system identification,telecommunication network
  routing,Telecommunication traffic,tomographic image
  reconstruction,tomography,Web and internet services,Web server}
  \strng{namehash}{CMHANR+1}
  \strng{fullhash}{CMHANRYB1}
  \field{abstract}{%
  Today's Internet is a massive, distributed network which continues to explode
  in size as e-commerce and related activities grow. The heterogeneous and
  largely unregulated structure of the Internet renders tasks such as dynamic
  routing, optimized service provision, service-level verification, and
  detection of anomalous/malicious behavior increasingly challenging tasks. The
  problem is compounded by the fact that one cannot rely on the cooperation of
  individual servers and routers to aid in the collection of network traffic
  measurements vital for these tasks. In many ways, network monitoring and
  inference problems bear a strong resemblance to other "inverse problems" in
  which key aspects of a system are not directly observable. Familiar signal
  processing problems such as tomographic image reconstruction, system
  identification, and array processing all have interesting interpretations in
  the networking context. This article introduces the new field of network
  tomography, a field which we believe will benefit greatly from the wealth of
  signal processing theory and algorithms%
  }
  \verb{doi}
  \verb 10.1109/79.998081
  \endverb
  \field{issn}{1053-5888}
  \field{number}{3}
  \field{pages}{47\bibrangedash 65}
  \field{title}{Internet Tomography}
  \field{volume}{19}
  \verb{file}
  \verb IEEE Xplore Abstract Record:D\:\\QUT\\Bibliography\\Zotero\\storage\\4A
  \verb A3DVNZ\\abstractReferences.html:;Coates et al_2002_Internet tomography.
  \verb pdf:D\:\\QUT\\Bibliography\\Zotero\\storage\\N4P6APQB\\Coates et al_200
  \verb 2_Internet tomography.pdf:application/pdf
  \endverb
  \field{journaltitle}{IEEE Signal Processing Magazine}
  \field{month}{05}
  \field{year}{2002}
\endentry

\entry{girard_performance_2006}{book}{}
  \name{author}{3}{}{%
    {{}%
     {Girard}{G.}%
     {André}{A.}%
     {}{}%
     {}{}}%
    {{}%
     {Sansò}{S.}%
     {Brunilde}{B.}%
     {}{}%
     {}{}}%
    {{}%
     {Vazquez-Abad}{V.-A.}%
     {Felida}{F.}%
     {}{}%
     {}{}}%
  }
  \list{publisher}{1}{%
    {{Springer}}%
  }
  \strng{namehash}{GASBVAF1}
  \strng{fullhash}{GASBVAF1}
  \field{title}{Performance Evaluation and Planning Methods for the next
  Generation Internet}
  \verb{url}
  \verb http://books.google.fr/books?hl=fr&lr=&id=Z_8f4NloZ0cC&oi=fnd&pg=PR11&d
  \verb q=Performance+Evaluation+and+Planning+Methods+for+the+Next+Generation+I
  \verb nternet&ots=tpg3jT4h5A&sig=zv2VPUfbEcjYwOhIocybctREF5o
  \endverb
  \field{volume}{6}
  \verb{file}
  \verb Girard et al_2006_Performance Evaluation and Planning Methods for the N
  \verb ext Generation Internet.pdf:D\:\\QUT\\Bibliography\\Zotero\\storage\\RW
  \verb JHTPA8\\Girard et al_2006_Performance Evaluation and Planning Methods f
  \verb or the Next Generation Internet.pdf:application/pdf;Snapshot:D\:\\QUT\\
  \verb Bibliography\\Zotero\\storage\\VUMMCADN\\books.html:
  \endverb
  \field{year}{2006}
  \field{urlday}{04}
  \field{urlmonth}{10}
  \field{urlyear}{2014}
\endentry

\entry{mardani_estimating_2015}{article}{}
  \name{author}{2}{}{%
    {{}%
     {Mardani}{M.}%
     {M.}{M.}%
     {}{}%
     {}{}}%
    {{}%
     {Giannakis}{G.}%
     {G.B.}{G.}%
     {}{}%
     {}{}}%
  }
  \keyw{Context,Convex optimization,Correlation,Estimation,IEEE
  transactions,low rank,nominal and anomalous traffic,Routing,Sparse
  matrices,Sparsity,spatiotemporal correlation,spatiotemporal phenomena}
  \strng{namehash}{MMGG1}
  \strng{fullhash}{MMGG1}
  \field{abstract}{%
  Mapping origin-destination (OD) network traffic is pivotal for network
  management and proactive security tasks. However, lack of sufficient
  flow-level measurements as well as potential anomalies pose major challenges
  towards this goal. Leveraging the spatiotemporal correlation of nominal
  traffic, and the sparse nature of anomalies, this paper brings forth a novel
  framework to map out nominal and anomalous traffic, which treats jointly
  important network monitoring tasks including traffic estimation, anomaly
  detection, and traffic interpolation. To this end, a convex program is first
  formulated with nuclear and -norm regularization to effect sparsity and low
  rank for the nominal and anomalous traffic with only the link counts and a
  small subset of OD-flow counts. Analysis and simulations confirm that the
  proposed estimator can exactly recover sufficiently low-dimensional nominal
  traffic and sporadic anomalies so long as the routing paths are sufficiently
  “spread-out” across the network, and an adequate amount of flow counts
  are randomly sampled. The results offer valuable insights about data
  acquisition strategies and network scenaria giving rise to accurate traffic
  estimation. For practical networks where the aforementioned conditions are
  possibly violated, the inherent spatiotemporal traffic patterns are taken
  into account by adopting a Bayesian approach along with a bilinear
  characterization of the nuclear and norms. The resultant nonconvex program
  involves quadratic regularizers with correlation matrices, learned
  systematically from (cyclo)stationary historical data.
  Alternating-minimization based algorithms with provable convergence are also
  developed to procure the estimates. Insightful tests with synthetic and real
  Internet data corroborate the effectiveness of the novel schemes.%
  }
  \verb{doi}
  \verb 10.1109/TNET.2015.2417809
  \endverb
  \field{issn}{1063-6692}
  \field{number}{99}
  \field{pages}{1\bibrangedash 15}
  \field{title}{Estimating Traffic and Anomaly Maps via Network Tomography}
  \verb{file}
  \verb IEEE Xplore Abstract Record:D\:\\QUT\\Bibliography\\Zotero\\storage\\KU
  \verb 6P4Z45\\freeabs_all.html:;Mardani and Giannakis - 2015 - Estimating Tra
  \verb ffic and Anomaly Maps via Network To.pdf:D\:\\QUT\\Bibliography\\Zotero
  \verb \\storage\\R6HWDZFR\\Mardani and Giannakis - 2015 - Estimating Traffic
  \verb and Anomaly Maps via Network To.pdf:application/pdf
  \endverb
  \field{journaltitle}{IEEE/ACM Transactions on Networking}
  \field{year}{2015}
\endentry

\entry{ortuzar_modelling_2011}{book}{}
  \name{author}{2}{}{%
    {{}%
     {Ortuzar}{O.}%
     {J.~D.}{J.~D.}%
     {}{}%
     {}{}}%
    {{}%
     {Willumsen}{W.}%
     {Luis~G.}{L.~G.}%
     {}{}%
     {}{}}%
  }
  \list{publisher}{1}{%
    {{John Wiley \& Sons, Ltd}}%
  }
  \strng{namehash}{OJDWLG1}
  \strng{fullhash}{OJDWLG1}
  \field{edition}{4th}
  \field{isbn}{978-0-470-76039-0}
  \field{pagetotal}{608}
  \field{title}{Modelling Transport}
  \verb{file}
  \verb Ortuzar_Willumsen_2011_Modelling Transport.pdf:D\:\\QUT\\Bibliography\\
  \verb Zotero\\storage\\DJEPR5QH\\Ortuzar_Willumsen_2011_Modelling Transport.p
  \verb df:application/pdf
  \endverb
  \field{year}{2011}
\endentry

\entry{wilson_entropy_1970}{book}{}
  \name{author}{1}{}{%
    {{}%
     {Wilson}{W.}%
     {A.~G.}{A.~G.}%
     {}{}%
     {}{}}%
  }
  \list{publisher}{1}{%
    {{Pion Ltd}}%
  }
  \strng{namehash}{WAG1}
  \strng{fullhash}{WAG1}
  \field{isbn}{0-85086-021-0}
  \field{pagetotal}{166}
  \field{title}{Entropy in Urban and Regional Modelling}
  \verb{url}
  \verb http://discovery.ucl.ac.uk/1315311/
  \endverb
  \verb{file}
  \verb Snapshot:D\:\\QUT\\Bibliography\\Zotero\\storage\\WGGQZRM9\\1315311.htm
  \verb l:
  \endverb
  \field{year}{1970}
  \field{urlday}{25}
  \field{urlmonth}{10}
  \field{urlyear}{2013}
\endentry

\entry{parry_estimation_2012}{article}{}
  \name{author}{2}{}{%
    {{}%
     {Parry}{P.}%
     {Katharina}{K.}%
     {}{}%
     {}{}}%
    {{}%
     {Hazelton}{H.}%
     {Martin~L.}{M.~L.}%
     {}{}%
     {}{}}%
  }
  \keyw{Boundary effects,Fisher information,Linear inverse problem,Maximum
  likelihood estimation,OD matrix,Poisson process}
  \strng{namehash}{PKHML1}
  \strng{fullhash}{PKHML1}
  \field{abstract}{%
  Estimation of origin–destination (OD) matrices from link count data is a
  challenging problem because of the highly indeterminate relationship between
  the observations and the latent route flows. Conversely, estimation is
  straightforward if we observe the path taken by each vehicle. We consider an
  intermediate problem of increasing practical importance, in which link count
  data is supplemented by routing information for a fraction of vehicles on the
  network. We develop a statistical model for these combined data sources and
  derive some tractable normal approximations thereof. We examine
  likelihood-based inference for these normal models under the assumption that
  the probability of vehicle tracking is known. We show that the likelihood
  theory can be non-standard because of boundary effects, and provide
  conditions under which such irregular behaviour will be observed in practice.
  For regular cases we outline connections with existing generalised least
  squares methods. We then consider estimation of OD matrices under estimated
  and/or misspecified models for the probability of vehicle tracking.
  Theoretical developments are complemented by simulation experiments and an
  illustrative example using a section of road network from the English city of
  Leicester.%
  }
  \verb{doi}
  \verb 10.1016/j.trb.2011.09.009
  \endverb
  \field{issn}{0191-2615}
  \field{number}{1}
  \field{pages}{175\bibrangedash 188}
  \field{shortjournal}{Transportation Research Part B: Methodological}
  \field{title}{Estimation of Origin–destination Matrices from Link Counts
  and Sporadic Routing Data}
  \field{volume}{46}
  \verb{file}
  \verb ScienceDirect Snapshot:D\:\\QUT\\Bibliography\\Zotero\\storage\\GI5DI6B
  \verb 4\\S019126151100138X.html:;Parry_Hazelton_2012_Estimation of origin–d
  \verb estination matrices from link counts and sporadic routing.pdf:D\:\\QUT\
  \verb \Bibliography\\Zotero\\storage\\W3JF7V8N\\Parry_Hazelton_2012_Estimatio
  \verb n of origin–destination matrices from link counts and sporadic routin
  \verb g.pdf:application/pdf
  \endverb
  \field{journaltitle}{Transportation Research Part B: Methodological}
  \field{month}{01}
  \field{year}{2012}
  \field{urlday}{12}
  \field{urlmonth}{03}
  \field{urlyear}{2014}
\endentry

\entry{kim_estimation_2010}{article}{}
  \name{author}{2}{}{%
    {{}%
     {Kim}{K.}%
     {Hyunmyung}{H.}%
     {}{}%
     {}{}}%
    {{}%
     {Jayakrishnan}{J.}%
     {R.}{R.}%
     {}{}%
     {}{}}%
  }
  \strng{namehash}{KHJR1}
  \strng{fullhash}{KHJR1}
  \field{abstract}{%
  Abstract In this paper we discuss a dynamic origin–destination (OD)
  estimation problem that has been used for identifying time-dependent travel
  demand on a road network. Even though a dynamic OD table is an indispensable
  data input for executing a dynamic traffic assignment, it is difficult to
  construct using the conventional OD construction method such as the four-step
  model. For this reason, a direct estimation method based on field traffic
  data such as link traffic counts has been used. However, the method does not
  account for a logical relationship between a travel demand pattern and
  socioeconomic attributes. In addition, the OD estimation method cannot
  guarantee the reliability of estimated results since the OD estimation
  problem has a property named the ‘underdetermined problem.’ In order to
  overcome such a problem, the method developed in this paper makes use of
  vehicle trajectory samples with link traffic counts. The new method is
  applied to numerical examples and shows promising capability for identifying
  a temporal and spatial travel demand pattern.%
  }
  \verb{doi}
  \verb 10.1080/03081060.2010.536629
  \endverb
  \field{issn}{0308-1060}
  \field{number}{8}
  \field{pages}{747\bibrangedash 768}
  \field{title}{The Estimation of a Time-Dependent {{OD}} Trip Table with
  Vehicle Trajectory Samples}
  \field{volume}{33}
  \verb{file}
  \verb Snapshot:D\:\\QUT\\Bibliography\\Zotero\\storage\\BBGNFUTF\\03081060.20
  \verb 10.html:;Kim_Jayakrishnan_2010_The estimation of a time-dependent OD tr
  \verb ip table with vehicle trajectory samples.pdf:D\:\\QUT\\Bibliography\\Zo
  \verb tero\\storage\\WP6R4XT6\\Kim_Jayakrishnan_2010_The estimation of a time
  \verb -dependent OD trip table with vehicle trajectory samples.pdf:applicatio
  \verb n/pdf
  \endverb
  \field{journaltitle}{Transportation Planning and Technology}
  \field{year}{2010}
\endentry

\entry{castillo_bayesian_2013}{article}{}
  \name{author}{4}{}{%
    {{}%
     {Castillo}{C.}%
     {Enrique}{E.}%
     {}{}%
     {}{}}%
    {{}%
     {Jiménez}{J.}%
     {Pilar}{P.}%
     {}{}%
     {}{}}%
    {{}%
     {Menéndez}{M.}%
     {José~María}{J.~M.}%
     {}{}%
     {}{}}%
    {{}%
     {Nogal}{N.}%
     {María}{M.}%
     {}{}%
     {}{}}%
  }
  \keyw{Conjugate priors,Economic Geography,Engineering Economics;
  Organization; Logistics; Marketing,Innovation/Technology
  Management,Origin–destination and link flow estimation,Prior assessment of
  hyperparameters,Regional/Spatial Science,Shifted-Gamma distribution}
  \strng{namehash}{CEJPMJM+1}
  \strng{fullhash}{CEJPMJMNM1}
  \field{abstract}{%
  In this paper a special conjugate Bayesian method, for reconstructing and
  estimating traffic flows, based on α-shifted-Gamma
  \ensuremath{\backslash}Upgamma(θ,λ) \ensuremath{\backslash}Upgamma
  (\ensuremath{\backslash}theta
  ,\ensuremath{\backslash},\ensuremath{\backslash}lambda ) models H(α,θ,λ)
  H(\ensuremath{\backslash}alpha
  ,\ensuremath{\backslash},\ensuremath{\backslash}theta
  ,\ensuremath{\backslash},\ensuremath{\backslash}lambda ) is given. If the
  numbers of users traveling through different routes are assumed to be
  independent H(α,θ,λ) H(\ensuremath{\backslash}alpha
  ,\ensuremath{\backslash},\ensuremath{\backslash}theta
  ,\ensuremath{\backslash},\ensuremath{\backslash}lambda) variables with common
  λ, \ensuremath{\backslash}lambda, the link, origin–destination (OD) and
  node flows are also H(α,θ,λ) H(\ensuremath{\backslash}alpha
  ,\ensuremath{\backslash},\ensuremath{\backslash}theta
  ,\ensuremath{\backslash},\ensuremath{\backslash}lambda ) random variables. We
  assume that the main source of information is plate scanning, which permits
  us to identify, totally or partially, the vehicle route, OD and link flows by
  scanning their corresponding plate numbers at an adequately selected subset
  of links. The reconstruction of the sample flows can be done exactly or
  approximately, depending on the intensity of the plate scanning sampling
  procedure. To this end a generalized least squares technique is used together
  with the conservation laws. A Bayesian approach using special conjugate
  families is proposed that allows us to estimate different traffic flows, such
  as route, OD-pair, scanned link or counted link flows. A detailed description
  of how the prior assessment, the sampling, the posterior updating and the
  obtention of the Bayesian distribution is given. Finally, one example of
  application is used to illustrate the methods and procedures.%
  }
  \verb{doi}
  \verb 10.1007/s11116-012-9443-4
  \endverb
  \field{issn}{0049-4488, 1572-9435}
  \field{number}{1}
  \field{pages}{173\bibrangedash 201}
  \field{shortjournal}{Transportation}
  \field{title}{A {{Bayesian}} Method for Estimating Traffic Flows Based on
  Plate Scanning}
  \field{volume}{40}
  \field{langid}{english}
  \verb{file}
  \verb Castillo et al_2013_A Bayesian method for estimating traffic flows base
  \verb d on plate scanning.pdf:D\:\\QUT\\Bibliography\\Zotero\\storage\\P92CR6
  \verb HS\\Castillo et al_2013_A Bayesian method for estimating traffic flows
  \verb based on plate scanning.pdf:application/pdf;Snapshot:D\:\\QUT\\Bibliogr
  \verb aphy\\Zotero\\storage\\TTGBE42P\\10.html:
  \endverb
  \field{journaltitle}{Transportation}
  \field{day}{01}
  \field{month}{01}
  \field{year}{2013}
\endentry

\entry{mellegard_origin/destination-estimation_2011}{inproceedings}{}
  \name{author}{3}{}{%
    {{}%
     {Mellegard}{M.}%
     {E.}{E.}%
     {}{}%
     {}{}}%
    {{}%
     {Moritz}{M.}%
     {S.}{S.}%
     {}{}%
     {}{}}%
    {{}%
     {Zahoor}{Z.}%
     {M.}{M.}%
     {}{}%
     {}{}}%
  }
  \keyw{Base stations,carriers networks,cellular network data,cellular
  radio,Clustering algorithms,Data privacy,end user
  privacy,geospatial,Hadoop,Land mobile radio cellular
  systems,location,mobility,MONOS
  devices,origin/destination-estimation,Origin-destination
  estimation,privacy,Roads,Smart phones}
  \strng{namehash}{MEMSZM1}
  \strng{fullhash}{MEMSZM1}
  \field{abstract}{%
  Today there are more than 600 billion geo special transactions every day in
  the US alone {[}1], and most of this data is passing through carriers
  networks. Hence, the carriers are sitting on a huge pile of potential
  knowledge which they could make more use of. Earlier attempt on this data
  have been very ambitious and often failed due to the nature of the data. In
  this paper present an innovative method that does enough, not too much with
  the data. Our method addresses how to handle the big data challenge as well
  as how one could plan to secure the privacy of the end users. The latter
  being one of the main reasons why carriers have not deployed anything similar
  yet. They have simply been too afraid of what would happen if the data would
  be mistreated or perceived as misused.%
  }
  \field{booktitle}{2011 {{IEEE}} 11th {{International Conference}} on {{Data
  Mining Workshops}} ({{ICDMW}})}
  \verb{doi}
  \verb 10.1109/ICDMW.2011.132
  \endverb
  \field{eventtitle}{2011 IEEE 11th International Conference on Data Mining
  Workshops (ICDMW)}
  \field{pages}{891\bibrangedash 896}
  \field{title}{Origin/destination-Estimation Using Cellular Network Data}
  \verb{file}
  \verb Mellegard et al_2011_Origin-Destination-estimation Using Cellular Netwo
  \verb rk Data.pdf:D\:\\QUT\\Bibliography\\Zotero\\storage\\8ZDU5JFB\\Mellegar
  \verb d et al_2011_Origin-Destination-estimation Using Cellular Network Data.
  \verb pdf:application/pdf;IEEE Xplore Abstract Record:D\:\\QUT\\Bibliography\
  \verb \Zotero\\storage\\J26M43DF\\abs_all.html:
  \endverb
  \field{month}{12}
  \field{year}{2011}
\endentry

\entry{iqbal_development_2014}{article}{}
  \name{author}{4}{}{%
    {{}%
     {Iqbal}{I.}%
     {Md.~Shahadat}{M.~S.}%
     {}{}%
     {}{}}%
    {{}%
     {Choudhury}{C.}%
     {Charisma~F.}{C.~F.}%
     {}{}%
     {}{}}%
    {{}%
     {Wang}{W.}%
     {Pu}{P.}%
     {}{}%
     {}{}}%
    {{}%
     {González}{G.}%
     {Marta~C.}{M.~C.}%
     {}{}%
     {}{}}%
  }
  \keyw{Mobile phone,Origin–destination,Traffic microsimulation,Video count}
  \strng{namehash}{IMSCCFWP+1}
  \strng{fullhash}{IMSCCFWPGMC1}
  \field{abstract}{%
  In this research, we propose a methodology to develop OD matrices using
  mobile phone Call Detail Records (CDR) and limited traffic counts. CDR, which
  consist of time stamped tower locations with caller IDs, are analyzed first
  and trips occurring within certain time windows are used to generate
  tower-to-tower transient OD matrices for different time periods. These are
  then associated with corresponding nodes of the traffic network and converted
  to node-to-node transient OD matrices. The actual OD matrices are derived by
  scaling up these node-to-node transient OD matrices. An optimization based
  approach, in conjunction with a microscopic traffic simulation platform, is
  used to determine the scaling factors that result best matches with the
  observed traffic counts. The methodology is demonstrated using CDR from 2.87
  million users of Dhaka, Bangladesh over a month and traffic counts from 13
  key locations over 3 days of that month. The applicability of the methodology
  is supported by a validation study.%
  }
  \verb{doi}
  \verb 10.1016/j.trc.2014.01.002
  \endverb
  \field{issn}{0968-090X}
  \field{pages}{63\bibrangedash 74}
  \field{shortjournal}{Transportation Research Part C: Emerging Technologies}
  \field{title}{Development of Origin-Destination Matrices Using Mobile Phone
  Call Data}
  \field{volume}{40}
  \verb{file}
  \verb Iqbal et al_2014_Development of origin–destination matrices using mob
  \verb ile phone call data.pdf:D\:\\QUT\\Bibliography\\Zotero\\storage\\PE6XDH
  \verb E2\\Iqbal et al_2014_Development of origin–destination matrices using
  \verb  mobile phone call data.pdf:application/pdf;ScienceDirect Snapshot:D\:\
  \verb \QUT\\Bibliography\\Zotero\\storage\\V835CWV2\\S0968090X14000059.html:
  \endverb
  \field{journaltitle}{Transportation Research Part C: Emerging Technologies}
  \field{month}{03}
  \field{year}{2014}
  \field{urlday}{23}
  \field{urlmonth}{02}
  \field{urlyear}{2014}
\endentry

\entry{alexander_origin-destination_2015}{article}{}
  \name{author}{4}{}{%
    {{}%
     {Alexander}{A.}%
     {Lauren}{L.}%
     {}{}%
     {}{}}%
    {{}%
     {Jiang}{J.}%
     {Shan}{S.}%
     {}{}%
     {}{}}%
    {{}%
     {Murga}{M.}%
     {Mikel}{M.}%
     {}{}%
     {}{}}%
    {{}%
     {González}{G.}%
     {Marta~C.}{M.~C.}%
     {}{}%
     {}{}}%
  }
  \keyw{Data mining,Human mobility,Mobile phone data,Travel surveys,Trip
  distribution,Trip production and attraction}
  \strng{namehash}{ALJSMM+1}
  \strng{fullhash}{ALJSMMGMC1}
  \field{abstract}{%
  In this work, we present methods to estimate average daily
  origin–destination trips from triangulated mobile phone records of millions
  of anonymized users. These records are first converted into clustered
  locations at which users engage in activities for an observed duration. These
  locations are inferred to be home, work, or other depending on observation
  frequency, day of week, and time of day, and represent a user’s origins and
  destinations. Since the arrival time and duration at these locations reflect
  the observed (based on phone usage) rather than true arrival time and
  duration of a user, we probabilistically infer departure time using survey
  data on trips in major US cities. Trips are then constructed for each user
  between two consecutive observations in a day. These trips are multiplied by
  expansion factors based on the population of a user’s home Census Tract and
  divided by the number of days on which we observed the user, distilling
  average daily trips. Aggregating individuals’ daily trips by Census Tract
  pair, hour of the day, and trip purpose results in trip matrices that form
  the basis for much of the analysis and modeling that inform transportation
  planning and investments. The applicability of the proposed methodology is
  supported by validation against the temporal and spatial distributions of
  trips reported in local and national surveys.%
  }
  \verb{doi}
  \verb 10.1016/j.trc.2015.02.018
  \endverb
  \field{issn}{0968-090X}
  \field{pages}{240\bibrangedash 250}
  \field{shortjournal}{Transportation Research Part C: Emerging Technologies}
  \field{title}{Origin-Destination Trips by Purpose and Time of Day Inferred
  from Mobile Phone Data}
  \field{volume}{58 part B}
  \verb{file}
  \verb Alexander et al. - 2015 - Origin–Destination Trips by Purpose and Tim
  \verb e of Da.pdf:D\:\\QUT\\Bibliography\\Zotero\\storage\\5RED8UEE\\Alexande
  \verb r et al. - 2015 - Origin–Destination Trips by Purpose and Time of Da.
  \verb pdf:application/pdf;ScienceDirect Snapshot:D\:\\QUT\\Bibliography\\Zote
  \verb ro\\storage\\FP3NH3NC\\S0968090X1500073X.html:
  \endverb
  \field{journaltitle}{Transportation Research Part C: Emerging Technologies}
  \field{month}{09}
  \field{year}{2015}
  \field{urlday}{25}
  \field{urlmonth}{03}
  \field{urlyear}{2015}
\endentry

\entry{herrera_evaluation_2010}{article}{}
  \name{author}{6}{}{%
    {{}%
     {Herrera}{H.}%
     {Juan~C.}{J.~C.}%
     {}{}%
     {}{}}%
    {{}%
     {Work}{W.}%
     {Daniel~B.}{D.~B.}%
     {}{}%
     {}{}}%
    {{}%
     {Herring}{H.}%
     {Ryan}{R.}%
     {}{}%
     {}{}}%
    {{}%
     {Ban}{B.}%
     {Xuegang~(Jeff)}{X.~J.}%
     {}{}%
     {}{}}%
    {{}%
     {Jacobson}{J.}%
     {Quinn}{Q.}%
     {}{}%
     {}{}}%
    {{}%
     {Bayen}{B.}%
     {Alexandre~M.}{A.~M.}%
     {}{}%
     {}{}}%
  }
  \keyw{GPS-enabled cell phones,Mobile sensors,Traffic monitoring
  systems,Traffic sensors}
  \strng{namehash}{HJCWDBHR+1}
  \strng{fullhash}{HJCWDBHRBXJJQBAM1}
  \field{abstract}{%
  The growing need of the driving public for accurate traffic information has
  spurred the deployment of large scale dedicated monitoring infrastructure
  systems, which mainly consist in the use of inductive loop detectors and
  video cameras. On-board electronic devices have been proposed as an
  alternative traffic sensing infrastructure, as they usually provide a
  cost-effective way to collect traffic data, leveraging existing communication
  infrastructure such as the cellular phone network. A traffic monitoring
  system based on GPS-enabled smartphones exploits the extensive coverage
  provided by the cellular network, the high accuracy in position and velocity
  measurements provided by GPS devices, and the existing infrastructure of the
  communication network. This article presents a field experiment nicknamed
  Mobile Century, which was conceived as a proof of concept of such a system.
  Mobile Century included 100 vehicles carrying a GPS-enabled Nokia N95 phone
  driving loops on a 10-mile stretch of I-880 near Union City, California, for
  8 h. Data were collected using virtual trip lines, which are geographical
  markers stored in the handset that probabilistically trigger position and
  speed updates when the handset crosses them. The proposed prototype system
  provided sufficient data for traffic monitoring purposes while managing the
  privacy of participants. The data obtained in the experiment were processed
  in real-time and successfully broadcast on the internet, demonstrating the
  feasibility of the proposed system for real-time traffic monitoring. Results
  suggest that a 2–3\% penetration of cell phones in the driver population is
  enough to provide accurate measurements of the velocity of the traffic flow.
  Data presented in this article can be downloaded from
  http://traffic.berkeley.edu.%
  }
  \verb{doi}
  \verb 10.1016/j.trc.2009.10.006
  \endverb
  \field{issn}{0968-090X}
  \field{number}{4}
  \field{pages}{568\bibrangedash 583}
  \field{shortjournal}{Transportation Research Part C: Emerging Technologies}
  \field{shorttitle}{Evaluation of Traffic Data Obtained via {{GPS}}-Enabled
  Mobile Phones}
  \field{title}{Evaluation of Traffic Data Obtained via {{GPS}}-Enabled Mobile
  Phones: The Mobile Century Field Experiment}
  \field{volume}{18}
  \verb{file}
  \verb ScienceDirect Snapshot:D\:\\QUT\\Bibliography\\Zotero\\storage\\GJ4UVKN
  \verb K\\S0968090X09001430.html:;Herrera et al_2010_Evaluation of traffic dat
  \verb a obtained via GPS-enabled mobile phones.pdf:D\:\\QUT\\Bibliography\\Zo
  \verb tero\\storage\\KRDMSP9Z\\Herrera et al_2010_Evaluation of traffic data
  \verb obtained via GPS-enabled mobile phones.pdf:application/pdf
  \endverb
  \field{journaltitle}{Transportation Research Part C: Emerging Technologies}
  \field{month}{08}
  \field{year}{2010}
  \field{urlday}{12}
  \field{urlmonth}{03}
  \field{urlyear}{2014}
\endentry

\entry{michau_bluetooth_2016}{article}{}
  \name{author}{6}{}{%
    {{}%
     {Michau}{M.}%
     {Gabriel}{G.}%
     {}{}%
     {}{}}%
    {{}%
     {Nantes}{N.}%
     {Alfredo}{A.}%
     {}{}%
     {}{}}%
    {{}%
     {Bhaskar}{B.}%
     {Ashish}{A.}%
     {}{}%
     {}{}}%
    {{}%
     {Chung}{C.}%
     {Edward}{E.}%
     {}{}%
     {}{}}%
    {{}%
     {Borgnat}{B.}%
     {Pierre}{P.}%
     {}{}%
     {}{}}%
    {{}%
     {Abry}{A.}%
     {Patrice}{P.}%
     {}{}%
     {}{}}%
  }
  \strng{namehash}{MGNABA+1}
  \strng{fullhash}{MGNABACEBPAP1}
  \field{title}{Bluetooth Data in Urban Context: Retrieving Vehicles
  Trajectories}
  \verb{file}
  \verb Michau et al_2016_Bluetooth Data in Urban Context.pdf:D\:\\QUT\\Bibliog
  \verb raphy\\Zotero\\storage\\KXD9JR4C\\Michau et al_2016_Bluetooth Data in U
  \verb rban Context.pdf:application/pdf
  \endverb
  \field{journaltitle}{submitted in IEEE Transaction on Intelligent Transport
  Systems}
  \field{year}{2016}
\endentry

\entry{hainen_estimating_2011}{article}{}
  \name{author}{6}{}{%
    {{}%
     {Hainen}{H.}%
     {Alexander}{A.}%
     {}{}%
     {}{}}%
    {{}%
     {Wasson}{W.}%
     {Jason}{J.}%
     {}{}%
     {}{}}%
    {{}%
     {Hubbard}{H.}%
     {Sarah}{S.}%
     {}{}%
     {}{}}%
    {{}%
     {Remias}{R.}%
     {Stephen}{S.}%
     {}{}%
     {}{}}%
    {{}%
     {Farnsworth}{F.}%
     {Grant}{G.}%
     {}{}%
     {}{}}%
    {{}%
     {Bullock}{B.}%
     {Darcy}{D.}%
     {}{}%
     {}{}}%
  }
  \strng{namehash}{HAWJHS+1}
  \strng{fullhash}{HAWJHSRSFGBD1}
  \field{abstract}{%
  Route choice is often assessed with either a modeling technique or field
  observations. Field observations have historically used a variation of
  license plate matching. The proposed technique assesses route choice and
  travel time that uses an anonymous Bluetooth media access control (MAC)
  address sampling technique as a surrogate for license plate matching to
  assess route choice. The Bluetooth sampling technique was used to evaluate
  the impact of an unexpected bridge closure in northwest Indiana, including an
  assessment of the proportion of vehicles using each of four alternate routes.
  The Bluetooth technology also provided a means to collect travel time data
  for each alternate route; these observed travel times were also compared with
  travel time estimates obtained by route classification and link distance. In
  general, the route choice behavior was consistent with observed travel time
  estimates. The Bluetooth sampling technique is cost-effective to deploy, and
  although results are approximate, direct measurement of travel times and
  route choice is useful for public agencies to assess mobility and travel time
  reliability along alternate routes.%
  }
  \verb{doi}
  \verb 10.3141/2256-06
  \endverb
  \field{pages}{43\bibrangedash 50}
  \field{title}{Estimating Route Choice and Travel Time Reliability with Field
  Observations of {{Bluetooth}} Probe Vehicles}
  \field{volume}{2256}
  \verb{file}
  \verb Hainen et al_2011_Estimating Route Choice and Travel Time Reliability w
  \verb ith Field Observations of.pdf:D\:\\QUT\\Bibliography\\Zotero\\storage\\
  \verb MXSU2KZ2\\Hainen et al_2011_Estimating Route Choice and Travel Time Rel
  \verb iability with Field Observations of.pdf:application/pdf;MetaPress Snaps
  \verb hot:D\:\\QUT\\Bibliography\\Zotero\\storage\\XUIA8UAM\\n257078847121068
  \verb .html:
  \endverb
  \field{journaltitle}{Transportation Research Record: Journal of the
  Transportation Research Board}
  \field{day}{01}
  \field{month}{12}
  \field{year}{2011}
  \field{urlday}{23}
  \field{urlmonth}{09}
  \field{urlyear}{2013}
\endentry

\entry{feng_vehicle_2015}{article}{}
  \name{author}{3}{}{%
    {{}%
     {Feng}{F.}%
     {Yu}{Y.}%
     {}{}%
     {}{}}%
    {{}%
     {Sun}{S.}%
     {Jian}{J.}%
     {}{}%
     {}{}}%
    {{}%
     {Chen}{C.}%
     {Peng}{P.}%
     {}{}%
     {}{}}%
  }
  \keyw{automatic vehicle identification,OD estimation,particle filter,traffic
  network,trajectory reconstruction}
  \strng{namehash}{FYSJCP1}
  \strng{fullhash}{FYSJCP1}
  \field{abstract}{%
  The origin–destination (OD) matrix and the vehicle trajectory data are
  critical to transportation planning, design, and operation management. On the
  basis of the deployment of automatic vehicle identification (AVI) technology
  in urban traffic networks in China, this study proposed a vehicle trajectory
  reconstruction method for a large-scale network by using AVI and traditional
  detector data. Particle filter theory was employed as the framework for this
  method that combines five spatial-temporal trajectory correction factors
  (i.e., the path consistency, the AVI measurability criterion, the travel time
  consistency, the gravity flow model, and the path-link flow matching model)
  to estimate the trajectory of a vehicle. The probabilities of the most likely
  trajectories are updated by the Bayesian method to approximate the ‘true’
  trajectory. The traffic network in the Beijing Olympic Park was selected as
  the test bed and was simulated by using VISSIM to create a complete set of
  vehicle trajectories. The accuracy of the resulting trajectory reconstruction
  exceeds 90\% when the AVI coverage is only 50\%, assuming an AVI detection
  error of 5\% for a closed network and an open network. The average relative
  error of a static OD matrix is 11.05\%. Although the accuracy of
  reconstruction exceeds 80\% when the AVI coverage is between 50\% and 40\%,
  the accuracy of a defective product-OD matrix decreases rapidly. The proposed
  method yields high estimation accuracy for the full trajectories of
  individual vehicles and the OD matrix, which demonstrates significant
  potential for traffic-related applications. Copyright © 2014 John Wiley \&
  Sons, Ltd.%
  }
  \verb{doi}
  \verb 10.1002/atr.1260
  \endverb
  \field{issn}{2042-3195}
  \field{number}{2}
  \field{pages}{174\bibrangedash 194}
  \field{shortjournal}{J. Adv. Transp.}
  \field{title}{Vehicle Trajectory Reconstruction Using Automatic Vehicle
  Identification and Traffic Count Data}
  \field{volume}{49}
  \field{langid}{english}
  \verb{file}
  \verb Feng et al_2015_Vehicle trajectory reconstruction using automatic vehic
  \verb le identification and.pdf:D\:\\QUT\\Bibliography\\Zotero\\storage\\PV3S
  \verb FA4D\\Feng et al_2015_Vehicle trajectory reconstruction using automatic
  \verb  vehicle identification and.pdf:application/pdf;Snapshot:D\:\\QUT\\Bibl
  \verb iography\\Zotero\\storage\\UAMAHR5M\\abstract.html:
  \endverb
  \field{journaltitle}{Journal of Advanced Transportation}
  \field{day}{01}
  \field{month}{03}
  \field{year}{2015}
  \field{urlday}{15}
  \field{urlmonth}{01}
  \field{urlyear}{2016}
\endentry

\entry{gomez_evaluation_2015}{article}{}
  \name{author}{3}{}{%
    {{}%
     {Gómez}{G.}%
     {Penélope}{P.}%
     {}{}%
     {}{}}%
    {{}%
     {Menéndez}{M.}%
     {Monica}{M.}%
     {}{}%
     {}{}}%
    {{}%
     {Mérida-Casermeiro}{M.-C.}%
     {Enrique}{E.}%
     {}{}%
     {}{}}%
  }
  \strng{namehash}{GPMMMCE1}
  \strng{fullhash}{GPMMMCE1}
  \field{abstract}{%
  In this paper, we evaluate the trade-offs between loop detector data and
  floating car data (FCD) for the real-time estimation of origin–destination
  (OD) matrices in small networks. The proposed methodology is based on a
  bi-level optimisation using fuzzy logic theory. Here we demonstrate that it
  provides accurate results with low computational cost, while presenting
  several advantages over other existing algorithms (especially in terms of
  data requirements, computational complexity, and quality of adjustment). The
  methodology is illustrated with three examples covering two different
  locations in the city of Zurich, Switzerland. Results are used to evaluate
  the trade-offs between loop detector coverage and the penetration rate of
  FCD, and to determine minimum values for ensuring a given accuracy level on
  the estimated OD matrices. In general, the resulting error in OD estimation
  is affected by the data redundancy in the network.%
  }
  \verb{doi}
  \verb 10.1080/21680566.2015.1025892
  \endverb
  \field{issn}{2168-0566}
  \field{number}{0}
  \field{pages}{1\bibrangedash 24}
  \field{title}{Evaluation of Trade-Offs between Two Data Sources for the
  Accurate Estimation of Origin-Destination Matrices}
  \field{volume}{0}
  \verb{file}
  \verb Snapshot:D\:\\QUT\\Bibliography\\Zotero\\storage\\2VBDJU7A\\21680566.20
  \verb 15.html:;Gómez et al_2015_Evaluation of trade-offs between two data so
  \verb urces for the accurate estimation.pdf:D\:\\QUT\\Bibliography\\Zotero\\s
  \verb torage\\7EA25KBG\\Gómez et al_2015_Evaluation of trade-offs between tw
  \verb o data sources for the accurate estimation.pdf:application/pdf
  \endverb
  \field{journaltitle}{Transportmetrica B: Transport Dynamics}
  \field{day}{26}
  \field{month}{03}
  \field{year}{2015}
  \field{urlday}{07}
  \field{urlmonth}{04}
  \field{urlyear}{2015}
\endentry

\entry{laharotte_spatiotemporal_2015}{article}{}
  \name{author}{6}{}{%
    {{}%
     {Laharotte}{L.}%
     {P.-A.}{P.-A.}%
     {}{}%
     {}{}}%
    {{}%
     {Billot}{B.}%
     {R.}{R.}%
     {}{}%
     {}{}}%
    {{}%
     {Come}{C.}%
     {E.}{E.}%
     {}{}%
     {}{}}%
    {{}%
     {Oukhellou}{O.}%
     {L.}{L.}%
     {}{}%
     {}{}}%
    {{}%
     {Nantes}{N.}%
     {A.}{A.}%
     {}{}%
     {}{}}%
    {{}%
     {El~Faouzi}{E.~F.}%
     {N.-E.}{N.-E.}%
     {}{}%
     {}{}}%
  }
  \keyw{$K$-means,arterial traffic,Australia,BFD classes,Bluetooth,Bluetooth
  data,Bluetooth detectors,Bluetooth fundamental diagram classes,Bluetooth
  origin-destination matrix,B-OD matrix,Brisbane urban area,clustering
  analysis,data processing,data reliability,data source,detection
  quality,Detectors,Estimation,forecasting,K-means algorithm,large-urban
  network,latent Dirichlet allocation,latent Dirichlet allocation
  (LDA),LDA,mode discrimination,motorway traffic,pattern clustering,penetration
  rate,road operators,road traffic capturing,road traffic
  control,spatiotemporal analysis,spatiotemporal clustering,spatiotemporal
  phenomena,spatiotemporal structure,temporal pattern detection,traffic
  condition monitoring,traffic conditions,traffic engineering computing,Traffic
  monitoring,traffic typology,travel time extraction,urban network,urban
  network characteristics,urban network dynamics,Vehicles,vehicle
  subpopulation}
  \strng{namehash}{LPABRCE+1}
  \strng{fullhash}{LPABRCEOLNAEFNE1}
  \field{abstract}{%
  The emergence of new technologies allows better monitoring of traffic
  conditions and understanding of urban network dynamics. Bluetooth technology
  is becoming widespread, as it represents a cost-effective means for capturing
  road traffic in both arterials and motorways. Although the extraction of
  travel time from Bluetooth data is fairly straightforward, data reliability
  and processing is still challenging with the issues of penetration rate, mode
  discrimination, and detection quality. This paper presents a methodological
  contribution to the use of Bluetooth data for the spatiotemporal analysis of
  a large urban network (Brisbane, Australia). It introduces the concept of the
  Bluetooth origin-destination (B-OD) matrix, which is built from a network of
  79 Bluetooth detectors located within the Brisbane urban area. The B-OD
  matrix describes the dynamics of a subpopulation of vehicles, between pairs
  of detectors. The results show that the characteristics of urban networks can
  be effectively represented through B-OD matrices. A comparison with loop
  detector data enables an assessment of the results' significance. Then, the
  spatiotemporal structure of the network is analyzed with two different
  clustering analyses, namely, latent Dirichlet allocation (LDA) and
  {\textdollar}K{\textdollar}-means. While LDA is used to detect a temporal
  pattern, the {\textdollar}K{\textdollar}-means algorithm highlights Bluetooth
  fundamental diagram (BFD) classes. The results show that Bluetooth data has
  the potential to be a reliable data source for traffic monitoring. By
  highlighting hidden structures of a large area, the algorithm outputs allow
  us to provide the road operators with a fine spatiotemporal analysis of their
  network, in terms of traffic conditions.%
  }
  \verb{doi}
  \verb 10.1109/TITS.2014.2367165
  \endverb
  \field{issn}{1524-9050}
  \field{number}{3}
  \field{pages}{1439\bibrangedash 1448}
  \field{shorttitle}{Spatiotemporal {{Analysis}} of {{Bluetooth Data}}}
  \field{title}{Spatiotemporal Analysis of {{Bluetooth}} Data: Application to a
  Large Urban Network}
  \field{volume}{16}
  \verb{file}
  \verb Laharotte et al_2015_Spatiotemporal Analysis of Bluetooth Data.pdf:D\:\
  \verb \QUT\\Bibliography\\Zotero\\storage\\3QXNRQ9I\\Laharotte et al_2015_Spa
  \verb tiotemporal Analysis of Bluetooth Data.pdf:application/pdf;IEEE Xplore
  \verb Abstract Record:D\:\\QUT\\Bibliography\\Zotero\\storage\\HSMWF9QP\\abs_
  \verb all.html:
  \endverb
  \field{journaltitle}{IEEE Transactions on Intelligent Transportation Systems}
  \field{month}{06}
  \field{year}{2015}
\endentry

\entry{combettes_douglas-rachford_2007}{article}{}
  \name{author}{2}{}{%
    {{}%
     {Combettes}{C.}%
     {Patrick~L.}{P.~L.}%
     {}{}%
     {}{}}%
    {{}%
     {Pesquet}{P.}%
     {Jean-Christophe}{J.-C.}%
     {}{}%
     {}{}}%
  }
  \strng{namehash}{CPLPJC1}
  \strng{fullhash}{CPLPJC1}
  \field{number}{4}
  \field{pages}{564\bibrangedash 574}
  \field{title}{A {{Douglas-Rachford}} Splitting Approach to Nonsmooth Convex
  Variational Signal Recovery}
  \verb{url}
  \verb http://ieeexplore.ieee.org/xpls/abs_all.jsp?arnumber=4407760
  \endverb
  \field{volume}{1}
  \verb{file}
  \verb Snapshot:D\:\\QUT\\Bibliography\\Zotero\\storage\\28WWXRTE\\abs_all.htm
  \verb l:;Combettes_Pesquet_2007_A Douglas–Rachford splitting approach to no
  \verb nsmooth convex variational signal.pdf:D\:\\QUT\\Bibliography\\Zotero\\s
  \verb torage\\5UJNANXX\\Combettes_Pesquet_2007_A Douglas–Rachford splitting
  \verb  approach to nonsmooth convex variational signal.pdf:application/pdf
  \endverb
  \field{journaltitle}{IEEE Journal of Selected Topics in Signal Processing}
  \field{year}{2007}
  \field{urlday}{20}
  \field{urlmonth}{03}
  \field{urlyear}{2015}
\endentry

\entry{chambolle_introduction_2010}{article}{}
  \name{author}{5}{}{%
    {{}%
     {Chambolle}{C.}%
     {Antonin}{A.}%
     {}{}%
     {}{}}%
    {{}%
     {Caselles}{C.}%
     {Vicent}{V.}%
     {}{}%
     {}{}}%
    {{}%
     {Cremers}{C.}%
     {Daniel}{D.}%
     {}{}%
     {}{}}%
    {{}%
     {Novaga}{N.}%
     {Matteo}{M.}%
     {}{}%
     {}{}}%
    {{}%
     {Pock}{P.}%
     {Thomas}{T.}%
     {}{}%
     {}{}}%
  }
  \strng{namehash}{CACVCD+1}
  \strng{fullhash}{CACVCDNMPT1}
  \field{pages}{263\bibrangedash 340}
  \field{title}{An Introduction to Total Variation for Image Analysis}
  \field{volume}{9}
  \verb{file}
  \verb Snapshot:D\:\\QUT\\Bibliography\\Zotero\\storage\\GE3SWBED\\books.html:
  \verb ;Chambolle et al_2010_An introduction to total variation for image anal
  \verb ysis.pdf:D\:\\QUT\\Bibliography\\Zotero\\storage\\ZBGUJB8T\\Chambolle e
  \verb t al_2010_An introduction to total variation for image analysis.pdf:app
  \verb lication/pdf
  \endverb
  \field{journaltitle}{Theoretical foundations and numerical methods for sparse
  recovery}
  \field{year}{2010}
\endentry

\entry{couprie_dual_2013}{article}{}
  \name{author}{5}{}{%
    {{}%
     {Couprie}{C.}%
     {Camille}{C.}%
     {}{}%
     {}{}}%
    {{}%
     {Grady}{G.}%
     {Leo}{L.}%
     {}{}%
     {}{}}%
    {{}%
     {Najman}{N.}%
     {Laurent}{L.}%
     {}{}%
     {}{}}%
    {{}%
     {Pesquet}{P.}%
     {Jean-Christophe}{J.-C.}%
     {}{}%
     {}{}}%
    {{}%
     {Talbot}{T.}%
     {Hugues}{H.}%
     {}{}%
     {}{}}%
  }
  \strng{namehash}{CCGLNL+1}
  \strng{fullhash}{CCGLNLPJCTH1}
  \field{number}{3}
  \field{pages}{1246\bibrangedash 1273}
  \field{title}{Dual Constrained {{TV}}-Based Regularization on Graphs}
  \field{volume}{6}
  \verb{file}
  \verb Couprie et al_2013_Dual constrained TV-based regularization on graphs.p
  \verb df:D\:\\QUT\\Bibliography\\Zotero\\storage\\4AZZJJUQ\\Couprie et al_201
  \verb 3_Dual constrained TV-based regularization on graphs.pdf:application/pd
  \verb f;Snapshot:D\:\\QUT\\Bibliography\\Zotero\\storage\\P2EM8C5H\\120895068
  \verb .html:
  \endverb
  \field{journaltitle}{SIAM Journal on Imaging Sciences}
  \field{year}{2013}
\endentry

\entry{pustelnik_wavelet-based_2016}{article}{}
  \name{author}{4}{}{%
    {{}%
     {Pustelnik}{P.}%
     {Nelly}{N.}%
     {}{}%
     {}{}}%
    {{}%
     {Benazza-Benhayia}{B.-B.}%
     {Amel}{A.}%
     {}{}%
     {}{}}%
    {{}%
     {Zheng}{Z.}%
     {Yuling}{Y.}%
     {}{}%
     {}{}}%
    {{}%
     {Pesquet}{P.}%
     {Jean-Christophe}{J.-C.}%
     {}{}%
     {}{}}%
  }
  \strng{namehash}{PNBBAZY+1}
  \strng{fullhash}{PNBBAZYPJC1}
  \field{pages}{1\bibrangedash 34}
  \field{title}{Wavelet-Based Image Deconvolution and Reconstruction}
  \verb{url}
  \verb https://hal.archives-ouvertes.fr/hal-01164833/
  \endverb
  \verb{file}
  \verb Snapshot:D\:\\QUT\\Bibliography\\Zotero\\storage\\7TXSTKGJ\\hal-0116483
  \verb 3.html:;Pustelnik et al. - 2016 - Wavelet-based Image Deconvolution and
  \verb  Reconstructi.pdf:D\:\\QUT\\Bibliography\\Zotero\\storage\\D47GKJHZ\\Pu
  \verb stelnik et al. - 2016 - Wavelet-based Image Deconvolution and Reconstru
  \verb cti.pdf:application/pdf
  \endverb
  \field{journaltitle}{Wiley Encyclopedia of Electrical and Electronics
  Engineering}
  \field{year}{2016}
  \field{urlday}{27}
  \field{urlmonth}{01}
  \field{urlyear}{2016}
\endentry

\entry{vu_splitting_2011}{article}{}
  \name{author}{1}{}{%
    {{}%
     {Vũ}{V.}%
     {Bằng~Công}{B.~C.}%
     {}{}%
     {}{}}%
  }
  \keyw{47H05,49M27,49M29,90C25,Algebra,Calculus of Variations and Optimal
  Control; Optimization,Cocoercivity,Composite
  operator,Duality,Forward-backward algorithm,Mathematics; general,Monotone
  inclusion,Monotone operator,Numeric Computing,Operator splitting,Primal-dual
  algorithm,Theory of Computation}
  \strng{namehash}{VBC1}
  \strng{fullhash}{VBC1}
  \field{abstract}{%
  We consider the problem of solving dual monotone inclusions involving sums of
  composite parallel-sum type operators. A feature of this work is to exploit
  explicitly the properties of the cocoercive operators appearing in the model.
  Several splitting algorithms recently proposed in the literature are
  recovered as special cases.%
  }
  \verb{doi}
  \verb 10.1007/s10444-011-9254-8
  \endverb
  \field{issn}{1019-7168, 1572-9044}
  \field{number}{3}
  \field{pages}{667\bibrangedash 681}
  \field{shortjournal}{Adv Comput Math}
  \field{title}{A Splitting Algorithm for Dual Monotone Inclusions Involving
  Cocoercive Operators}
  \field{volume}{38}
  \field{langid}{english}
  \verb{file}
  \verb Snapshot:D\:\\QUT\\Bibliography\\Zotero\\storage\\3N7N89BC\\s10444-011-
  \verb 9254-8.html:;Vũ_2011_A splitting algorithm for dual monotone inclusion
  \verb s involving cocoercive.pdf:D\:\\QUT\\Bibliography\\Zotero\\storage\\9CN
  \verb H3WBG\\Vũ_2011_A splitting algorithm for dual monotone inclusions invo
  \verb lving cocoercive.pdf:application/pdf
  \endverb
  \field{journaltitle}{Advances in Computational Mathematics}
  \field{day}{29}
  \field{month}{11}
  \field{year}{2011}
  \field{urlday}{27}
  \field{urlmonth}{01}
  \field{urlyear}{2016}
\endentry

\entry{combettes_primal-dual_2011}{article}{}
  \name{author}{2}{}{%
    {{}%
     {Combettes}{C.}%
     {Patrick~L.}{P.~L.}%
     {}{}%
     {}{}}%
    {{}%
     {Pesquet}{P.}%
     {Jean-Christophe}{J.-C.}%
     {}{}%
     {}{}}%
  }
  \keyw{47H05,49M27,49M29,49N15,90C25,Analysis,Geometry,Maximal monotone
  operator,Monotone inclusion,Nonsmooth convex optimization,Parallel
  sum,Set-valued duality,Splitting algorithm}
  \strng{namehash}{CPLPJC1}
  \strng{fullhash}{CPLPJC1}
  \field{abstract}{%
  We propose a primal-dual splitting algorithm for solving monotone inclusions
  involving a mixture of sums, linear compositions, and parallel sums of
  set-valued and Lipschitzian operators. An important feature of the algorithm
  is that the Lipschitzian operators present in the formulation can be
  processed individually via explicit steps, while the set-valued operators are
  processed individually via their resolvents. In addition, the algorithm is
  highly parallel in that most of its steps can be executed simultaneously.
  This work brings together and notably extends various types of structured
  monotone inclusion problems and their solution methods. The application to
  convex minimization problems is given special attention.%
  }
  \verb{doi}
  \verb 10.1007/s11228-011-0191-y
  \endverb
  \field{issn}{1877-0533, 1877-0541}
  \field{number}{2}
  \field{pages}{307\bibrangedash 330}
  \field{shortjournal}{Set-Valued Anal}
  \field{title}{Primal-Dual Splitting Algorithm for Solving Inclusions with
  Mixtures of Composite, {{Lipschitzian}}, and Parallel-Sum Type Monotone
  Operators}
  \field{volume}{20}
  \field{langid}{english}
  \verb{file}
  \verb Snapshot:D\:\\QUT\\Bibliography\\Zotero\\storage\\U6RZ2IG5\\s11228-011-
  \verb 0191-y.html:;Combettes_Pesquet_2011_Primal-Dual Splitting Algorithm for
  \verb  Solving Inclusions with Mixtures of.pdf:D\:\\QUT\\Bibliography\\Zotero
  \verb \\storage\\ZSZDCNWU\\Combettes_Pesquet_2011_Primal-Dual Splitting Algor
  \verb ithm for Solving Inclusions with Mixtures of.pdf:application/pdf
  \endverb
  \field{journaltitle}{Set-Valued and Variational Analysis}
  \field{day}{27}
  \field{month}{08}
  \field{year}{2011}
  \field{urlday}{27}
  \field{urlmonth}{01}
  \field{urlyear}{2016}
\endentry

\entry{condat_primal-dual_2013}{article}{}
  \name{author}{1}{}{%
    {{}%
     {Condat}{C.}%
     {Laurent}{L.}%
     {}{}%
     {}{}}%
  }
  \keyw{Applications of Mathematics,Calculus of Variations and Optimal Control;
  Optimization,Convex and nonsmooth optimization,Douglas–Rachford
  method,Engineering; general,Fenchel–Rockafellar duality,Forward–backward
  method,Monotone inclusion,Operations Research/Decision Theory,Operator
  splitting,Optimization,Primal–dual algorithm,Proximal method,Theory of
  Computation}
  \strng{namehash}{CL1}
  \strng{fullhash}{CL1}
  \field{abstract}{%
  We propose a new first-order splitting algorithm for solving jointly the
  primal and dual formulations of large-scale convex minimization problems
  involving the sum of a smooth function with Lipschitzian gradient, a
  nonsmooth proximable function, and linear composite functions. This is a full
  splitting approach, in the sense that the gradient and the linear operators
  involved are applied explicitly without any inversion, while the nonsmooth
  functions are processed individually via their proximity operators. This work
  brings together and notably extends several classical splitting schemes, like
  the forward–backward and Douglas–Rachford methods, as well as the recent
  primal–dual method of Chambolle and Pock designed for problems with linear
  composite terms.%
  }
  \verb{doi}
  \verb 10.1007/s10957-012-0245-9
  \endverb
  \field{issn}{0022-3239, 1573-2878}
  \field{number}{2}
  \field{pages}{460\bibrangedash 479}
  \field{shortjournal}{J Optim Theory Appl}
  \field{title}{A Primal-Dual Splitting Method for Convex Optimization
  Involving {{Lipschitzian}}, Proximable and Linear Composite Terms}
  \field{volume}{158}
  \field{langid}{english}
  \verb{file}
  \verb Snapshot:D\:\\QUT\\Bibliography\\Zotero\\storage\\5F5J4NWA\\10.html:;Co
  \verb ndat_2013_A Primal–Dual Splitting Method for Convex Optimization Invo
  \verb lving Lipschitzian,.pdf:D\:\\QUT\\Bibliography\\Zotero\\storage\\TT869I
  \verb ME\\Condat_2013_A Primal–Dual Splitting Method for Convex Optimizatio
  \verb n Involving Lipschitzian,.pdf:application/pdf
  \endverb
  \field{journaltitle}{Journal of Optimization Theory and Applications}
  \field{day}{01}
  \field{month}{08}
  \field{year}{2013}
  \field{urlday}{10}
  \field{urlmonth}{01}
  \field{urlyear}{2014}
\endentry

\entry{komodakis_playing_2015}{article}{}
  \name{author}{2}{}{%
    {{}%
     {Komodakis}{K.}%
     {N.}{N.}%
     {}{}%
     {}{}}%
    {{}%
     {Pesquet}{P.}%
     {J.-C.}{J.-C.}%
     {}{}%
     {}{}}%
  }
  \keyw{Algorithm design and analysis,computer vision,Context modeling,convex
  analysis,Convex functions,discrete optimization,duality
  (mathematics),high-dimensional problems,image processing,large-scale
  optimization problem solving,machine learning,nonsmooth
  optimization,numerical analysis,numerical
  methods,optimisation,Optimization,parallel algorithms,parallel
  processing,primal-dual approaches,signal-image processing,Signal processing
  algorithms}
  \strng{namehash}{KNPJC1}
  \strng{fullhash}{KNPJC1}
  \field{abstract}{%
  Optimization methods are at the core of many problems in signal/image
  processing, computer vision, and machine learning. For a long time, it has
  been recognized that looking at the dual of an optimization problem may
  drastically simplify its solution. However, deriving efficient strategies
  that jointly bring into play the primal and dual problems is a more recent
  idea that has generated many important new contributions in recent years.
  These novel developments are grounded in the recent advances in convex
  analysis, discrete optimization, parallel processing, and nonsmooth
  optimization with an emphasis on sparsity issues. In this article, we aim to
  present the principles of primal-dual approaches while providing an overview
  of the numerical methods that have been proposed in different contexts. Last
  but not least, primal-dual methods lead to algorithms that are easily
  parallelizable. Today, such parallel algorithms are becoming increasingly
  important for efficiently handling high-dimensional problems.%
  }
  \verb{doi}
  \verb 10.1109/MSP.2014.2377273
  \endverb
  \field{issn}{1053-5888}
  \field{number}{6}
  \field{pages}{31\bibrangedash 54}
  \field{shorttitle}{Playing with {{Duality}}}
  \field{title}{Playing with Duality: An Overview of Recent Primal-Dual
  Approaches for Solving Large-Scale Optimization Problems}
  \field{volume}{32}
  \verb{file}
  \verb Komodakis_Pesquet_2015_Playing with Duality.pdf:D\:\\QUT\\Bibliography\
  \verb \Zotero\\storage\\9WVBVDS6\\Komodakis_Pesquet_2015_Playing with Duality
  \verb .pdf:application/pdf;IEEE Xplore Abstract Record:D\:\\QUT\\Bibliography
  \verb \\Zotero\\storage\\QBVITGCP\\articleDetails.html:
  \endverb
  \field{journaltitle}{IEEE Signal Processing Magazine}
  \field{month}{11}
  \field{year}{2015}
\endentry

\entry{moreau_proximite_1965}{article}{}
  \name{author}{1}{}{%
    {{}%
     {Moreau}{M.}%
     {Jean-Jacques}{J.-J.}%
     {}{}%
     {}{}}%
  }
  \strng{namehash}{MJJ1}
  \strng{fullhash}{MJJ1}
  \field{pages}{273\bibrangedash 299}
  \field{title}{Proximité et Dualité Dans Un Espace {{Hilbertien}}}
  \verb{url}
  \verb http://archive.numdam.org/article/BSMF_1965__93__273_0.pdf
  \endverb
  \field{volume}{93}
  \verb{file}
  \verb Moreau_1965_Proximité et dualité dans un espace hilbertien.pdf:D\:\\Q
  \verb UT\\Bibliography\\Zotero\\storage\\3WHQDE34\\Moreau_1965_Proximité et
  \verb dualité dans un espace hilbertien.pdf:application/pdf
  \endverb
  \field{journaltitle}{Bulletin de la Société mathématique de France}
  \field{year}{1965}
  \field{urlday}{21}
  \field{urlmonth}{01}
  \field{urlyear}{2016}
\endentry

\entry{combettes_proximal_2010}{incollection}{}
  \name{author}{2}{}{%
    {{}%
     {Combettes}{C.}%
     {Patrick~L.}{P.~L.}%
     {}{}%
     {}{}}%
    {{}%
     {Pesquet}{P.}%
     {Jean~Christophe}{J.~C.}%
     {}{}%
     {}{}}%
  }
  \name{editor}{6}{}{%
    {{}%
     {Bauschke}{B.}%
     {Heinz~H.}{H.~H.}%
     {}{}%
     {}{}}%
    {{}%
     {Burachik}{B.}%
     {R.}{R.}%
     {}{}%
     {}{}}%
    {{}%
     {Combettes}{C.}%
     {Patrick~L.}{P.~L.}%
     {}{}%
     {}{}}%
    {{}%
     {Elser}{E.}%
     {V.}{V.}%
     {}{}%
     {}{}}%
    {{}%
     {Luke}{L.}%
     {D.~R.}{D.~R.}%
     {}{}%
     {}{}}%
    {{}%
     {Wolkowicz}{W.}%
     {H.}{H.}%
     {}{}%
     {}{}}%
  }
  \list{publisher}{1}{%
    {{Springer-Verlag}}%
  }
  \strng{namehash}{CPLPJC1}
  \strng{fullhash}{CPLPJC1}
  \field{booktitle}{Fixed-{{Point Algorithms}} for {{Inverse Problems}} in
  {{Science}} and {{Engineering}}}
  \field{pages}{185\bibrangedash 212}
  \field{title}{Proximal Splitting Methods in Signal Processing}
  \verb{url}
  \verb http://link.springer.com/chapter/10.1007/978-1-4419-9569-8_10
  \endverb
  \list{location}{1}{%
    {{New-York}}%
  }
  \verb{file}
  \verb Combettes_Pesquet_2011_Proximal splitting methods in signal processing.
  \verb pdf:D\:\\QUT\\Bibliography\\Zotero\\storage\\2XSUS64T\\Combettes_Pesque
  \verb t_2011_Proximal splitting methods in signal processing.pdf:application/
  \verb pdf;Snapshot:D\:\\QUT\\Bibliography\\Zotero\\storage\\8MACQC56\\978-1-4
  \verb 419-9569-8_10.html:
  \endverb
  \field{year}{2010}
  \field{urlday}{29}
  \field{urlmonth}{09}
  \field{urlyear}{2014}
\endentry

\entry{bauschke_convex_2011}{book}{}
  \name{author}{2}{}{%
    {{}%
     {Bauschke}{B.}%
     {Heinz~H.}{H.~H.}%
     {}{}%
     {}{}}%
    {{}%
     {Combettes}{C.}%
     {Patrick~L.}{P.~L.}%
     {}{}%
     {}{}}%
  }
  \list{publisher}{1}{%
    {{Springer}}%
  }
  \strng{namehash}{BHHCPL1}
  \strng{fullhash}{BHHCPL1}
  \field{isbn}{978-1-4419-9466-0}
  \field{series}{CMS Books in Mathematics Ser.}
  \field{title}{Convex Analysis and Monotone Operator Theory in {{Hilbert}}
  Spaces}
  \verb{file}
  \verb convex analysis and monotone operator theory in hilbert spaces - QUT Li
  \verb brary Quick Find:D\:\\QUT\\Bibliography\\Zotero\\storage\\26SMSQVM\\sea
  \verb rch.html:;Bauschke_Combettes_2011_Convex Analysis and Monotone Operator
  \verb  Theory in Hilbert Spaces.pdf:D\:\\QUT\\Bibliography\\Zotero\\storage\\
  \verb 3WT56ZBD\\Bauschke_Combettes_2011_Convex Analysis and Monotone Operator
  \verb  Theory in Hilbert Spaces.pdf:application/pdf
  \endverb
  \field{year}{2011}
\endentry

\entry{parikh_proximal_2013}{article}{}
  \name{author}{2}{}{%
    {{}%
     {Parikh}{P.}%
     {Neal}{N.}%
     {}{}%
     {}{}}%
    {{}%
     {Boyd}{B.}%
     {Stephen}{S.}%
     {}{}%
     {}{}}%
  }
  \strng{namehash}{PNBS1}
  \strng{fullhash}{PNBS1}
  \field{number}{3}
  \field{pages}{123\bibrangedash 231}
  \field{title}{Proximal Algorithms}
  \verb{url}
  \verb http://web.stanford.edu/~boyd/papers/pdf/prox_algs.pdf
  \endverb
  \field{volume}{1}
  \verb{file}
  \verb Parikh_Boyd_2013_Proximal Algorithms.pdf:D\:\\QUT\\Bibliography\\Zotero
  \verb \\storage\\XXBPBFXR\\Parikh_Boyd_2013_Proximal Algorithms.pdf:applicati
  \verb on/pdf
  \endverb
  \field{journaltitle}{Foundations and Trends in optimization}
  \field{year}{2013}
  \field{urlday}{20}
  \field{urlmonth}{03}
  \field{urlyear}{2015}
\endentry

\entry{theodoridis_adaptive_2011}{article}{}
  \name{author}{3}{}{%
    {{}%
     {Theodoridis}{T.}%
     {S.}{S.}%
     {}{}%
     {}{}}%
    {{}%
     {Slavakis}{S.}%
     {K.}{K.}%
     {}{}%
     {}{}}%
    {{}%
     {Yamada}{Y.}%
     {I.}{I.}%
     {}{}%
     {}{}}%
  }
  \keyw{adaptive equalization,adaptive learning,adaptive robust
  beamforming,Array signal processing,Estimation,function estimation,Hilbert
  spaces,Learning systems,Least squares approximation,numerical
  simulations,Optimization,parameter estimation,POCS,projections onto convex
  sets,reproducing kernel Hilbert spaces,RKHS,set theory,signal
  processing,Signal processing algorithms,Training data}
  \strng{namehash}{TSSKYI1}
  \strng{fullhash}{TSSKYI1}
  \field{abstract}{%
  This article presents a general tool for convexly constrained
  parameter/function estimation both for classification and regression tasks,
  in a timeadaptive setting and in (infinite dimensional) reproducing kernel
  Hilbert spaces (RKHS). The thematical framework is that of the set theoretic
  estimation formulation and the classical projections onto convex sets (POCS)
  theory. However, in contrast to the classical POCS methodology, which assumes
  a finite number of convex sets, our method builds upon our recent extension
  of the theory, which considers an infinite number of convex sets. Such a
  context is necessary to cope with the adaptive setting rationale, where data
  arrive sequentially. This article's goal is to review the advances that have
  taken place in this area over the years and present them, in simple geometric
  arguments, as an integral part and natural evolution of the classical POCS
  methodology. The structure of the resulting algorithms is such that it allows
  extension to general RKHS. In this perspective, two very powerful techniques,
  convex optimization and (implicit) mapping to RKHS, are combined, which
  provide a framework for a unifying treatment of linear and nonlinear modeling
  of both classification and regression tasks. Typical signal processing
  problems, such as filtering, smoothing, equalization, and beamforming, fall
  under this common umbrella. The methodology allows for the incorporation of a
  set of convex constraints, which encode a priori information. Convexity,
  rather than differentiability, is the only prerequisite for adopting error
  measures that quantify the model's fit against a set of training data points.
  Moreover, the complexity per iteration step remains linear with respect to
  the number of unknown parameters. The potential of the theory is demonstrated
  via numerical simulations for two typical problems; adaptive equalization and
  adaptive robust beamforming.%
  }
  \verb{doi}
  \verb 10.1109/MSP.2010.938752
  \endverb
  \field{issn}{1053-5888}
  \field{number}{1}
  \field{pages}{97\bibrangedash 123}
  \field{title}{Adaptive Learning in a World of Projections}
  \field{volume}{28}
  \verb{file}
  \verb IEEE Xplore Abstract Record:D\:\\QUT\\Bibliography\\Zotero\\storage\\D8
  \verb 555SRQ\\articleDetails.html:;Theodoridis et al_2011_Adaptive Learning i
  \verb n a World of Projections.pdf:D\:\\QUT\\Bibliography\\Zotero\\storage\\F
  \verb ZF9IBSN\\Theodoridis et al_2011_Adaptive Learning in a World of Project
  \verb ions.pdf:application/pdf
  \endverb
  \field{journaltitle}{IEEE Signal Processing Magazine}
  \field{year}{2011}
\endentry

\entry{chaux_nested_2009}{article}{}
  \name{author}{3}{}{%
    {{}%
     {Chaux}{C.}%
     {Caroline}{C.}%
     {}{}%
     {}{}}%
    {{}%
     {Pesquet}{P.}%
     {Jean-Christophe}{J.-C.}%
     {}{}%
     {}{}}%
    {{}%
     {Pustelnik}{P.}%
     {Nelly}{N.}%
     {}{}%
     {}{}}%
  }
  \keyw{Engineering--Mechanical Engineering,Physics--Optics}
  \strng{namehash}{CCPJCPN1}
  \strng{fullhash}{CCPJCPN1}
  \field{abstract}{%
  The objective of this paper is to develop methods for solving image recovery
  problems subject to constraints on the solution. More precisely, we will be
  interested in problems which can be formulated as the minimization over a
  closed convex constraint set of the sum of two convex functions
  {\textdollar}f{\textdollar} and {\textdollar}g{\textdollar}, where
  {\textdollar}f{\textdollar} may be nonsmooth and {\textdollar}g{\textdollar}
  is differentiable with a Lipschitz-continuous gradient. To reach this goal,
  we derive two types of algorithms that combine forward-backward and
  Douglas-Rachford iterations. The weak convergence of the proposed algorithms
  is proved. In the case when the Lipschitz-continuity property of the gradient
  of {\textdollar}g{\textdollar} is not satisfied, we also show that, under
  some assumptions, it remains possible to apply these methods to the
  considered optimization problem by making use of a quadratic extension
  technique. The effectiveness of the algorithms is demonstrated for two
  wavelet-based image restoration problems involving a signal-dependent
  Gaussian noise and a Poisson noise, respectively.%
  }
  \verb{doi}
  \verb http://dx.doi.org.ezp01.library.qut.edu.au/10.1137/080727749
  \endverb
  \field{issn}{19364954}
  \field{number}{2}
  \field{pages}{33}
  \field{title}{Nested Iterative Algorithms for Convex Constrained Image
  Recovery Problems}
  \field{volume}{2}
  \field{langid}{english}
  \verb{file}
  \verb Chaux et al_2009_Nested Iterative Algorithms for Convex Constrained Ima
  \verb ge Recovery Problems.pdf:D\:\\QUT\\Bibliography\\Zotero\\storage\\2RRFI
  \verb 76Z\\Chaux et al_2009_Nested Iterative Algorithms for Convex Constraine
  \verb d Image Recovery Problems.pdf:application/pdf;Snapshot:D\:\\QUT\\Biblio
  \verb graphy\\Zotero\\storage\\BCU6BH28\\1.html:
  \endverb
  \field{journaltitle}{SIAM Journal on Imaging Sciences}
  \field{year}{2009}
  \field{urlday}{20}
  \field{urlmonth}{03}
  \field{urlyear}{2015}
\endentry

\entry{chaux_variational_2007}{article}{}
  \name{author}{4}{}{%
    {{}%
     {Chaux}{C.}%
     {Caroline}{C.}%
     {}{}%
     {}{}}%
    {{}%
     {Combettes}{C.}%
     {Patrick~L.}{P.~L.}%
     {}{}%
     {}{}}%
    {{}%
     {Pesquet}{P.}%
     {Jean-Christophe}{J.-C.}%
     {}{}%
     {}{}}%
    {{}%
     {Wajs}{W.}%
     {Valérie~R.}{V.~R.}%
     {}{}%
     {}{}}%
  }
  \strng{namehash}{CCCPLPJC+1}
  \strng{fullhash}{CCCPLPJCWVR1}
  \verb{doi}
  \verb 10.1088/0266-5611/23/4/008
  \endverb
  \field{number}{4}
  \field{pages}{1495\bibrangedash 1518}
  \field{title}{A Variational Formulation for Frame-Based Inverse Problems}
  \field{volume}{23}
  \verb{file}
  \verb Chaux et al_2007_A variational formulation for frame-based inverse prob
  \verb lems.pdf:D\:\\QUT\\Bibliography\\Zotero\\storage\\9CIFCP9K\\Chaux et al
  \verb _2007_A variational formulation for frame-based inverse problems.pdf:ap
  \verb plication/pdf;Snapshot:D\:\\QUT\\Bibliography\\Zotero\\storage\\KSTX5FD
  \verb F\\meta\;jsessionid=2CA4550902CC4FEF2C41C9DF1448A0F8.c5.iopscience.cld.
  \verb iop.html:
  \endverb
  \field{journaltitle}{Inverse Problems}
  \field{year}{2007}
  \field{urlday}{21}
  \field{urlmonth}{01}
  \field{urlyear}{2016}
\endentry

\entry{pustelnik_parallel_2009}{article}{}
  \name{author}{3}{}{%
    {{}%
     {Pustelnik}{P.}%
     {Nelly}{N.}%
     {}{}%
     {}{}}%
    {{}%
     {Chaux}{C.}%
     {Caroline}{C.}%
     {}{}%
     {}{}}%
    {{}%
     {Pesquet}{P.}%
     {Jean-Christophe}{J.-C.}%
     {}{}%
     {}{}}%
  }
  \strng{namehash}{PNCCPJC1}
  \strng{fullhash}{PNCCPJC1}
  \field{title}{Parallel Proximal Algorithm for Image Restoration Using Hybrid
  Regularization -- Extended Version}
  \verb{url}
  \verb http://arxiv.org/abs/0911.1536
  \endverb
  \verb{file}
  \verb Snapshot:D\:\\QUT\\Bibliography\\Zotero\\storage\\ADAFXBEG\\0911.html:;
  \verb Pustelnik et al. - 2009 - Parallel Proximal Algorithm for Image Restora
  \verb tion .pdf:D\:\\QUT\\Bibliography\\Zotero\\storage\\ZWNZW2UU\\Pustelnik
  \verb et al. - 2009 - Parallel Proximal Algorithm for Image Restoration .pdf:
  \verb application/pdf
  \endverb
  \field{day}{08}
  \field{month}{11}
  \field{year}{2009}
  \field{urlday}{18}
  \field{urlmonth}{11}
  \field{urlyear}{2015}
\endentry

\entry{donoho_-noising_1995}{article}{}
  \name{author}{1}{}{%
    {{}%
     {Donoho}{D.}%
     {David~L.}{D.~L.}%
     {}{}%
     {}{}}%
  }
  \strng{namehash}{DDL1}
  \strng{fullhash}{DDL1}
  \field{number}{3}
  \field{pages}{613\bibrangedash 627}
  \field{title}{De-Noising by Soft-Thresholding}
  \verb{url}
  \verb http://ieeexplore.ieee.org/xpls/abs_all.jsp?arnumber=382009
  \endverb
  \field{volume}{41}
  \verb{file}
  \verb Snapshot:D\:\\QUT\\Bibliography\\Zotero\\storage\\78HANRHB\\freeabs_all
  \verb .html:;Donoho_1995_De-noising by soft-thresholding.pdf:D\:\\QUT\\Biblio
  \verb graphy\\Zotero\\storage\\BD5R2KMH\\Donoho_1995_De-noising by soft-thres
  \verb holding.pdf:application/pdf
  \endverb
  \field{journaltitle}{Information Theory, IEEE Transactions on}
  \field{year}{1995}
  \field{urlday}{27}
  \field{urlmonth}{01}
  \field{urlyear}{2016}
\endentry

\entry{combettes_signal_2005}{article}{}
  \name{author}{2}{}{%
    {{}%
     {Combettes}{C.}%
     {Patrick~L.}{P.~L.}%
     {}{}%
     {}{}}%
    {{}%
     {Wajs}{W.}%
     {Valérie~R.}{V.~R.}%
     {}{}%
     {}{}}%
  }
  \strng{namehash}{CPLWVR1}
  \strng{fullhash}{CPLWVR1}
  \field{number}{4}
  \field{pages}{1168\bibrangedash 1200}
  \field{title}{Signal Recovery by Proximal Forward-Backward Splitting}
  \verb{url}
  \verb http://epubs.siam.org/doi/abs/10.1137/050626090
  \endverb
  \field{volume}{4}
  \verb{file}
  \verb Combettes_Wajs_2005_Signal recovery by proximal forward-backward splitt
  \verb ing.pdf:D\:\\QUT\\Bibliography\\Zotero\\storage\\5PJXK5B8\\Combettes_Wa
  \verb js_2005_Signal recovery by proximal forward-backward splitting.pdf:appl
  \verb ication/pdf;Snapshot:D\:\\QUT\\Bibliography\\Zotero\\storage\\VI9WRQJQ\
  \verb \050626090.html:
  \endverb
  \field{journaltitle}{Multiscale Modeling \& Simulation}
  \field{year}{2005}
  \field{urlday}{29}
  \field{urlmonth}{09}
  \field{urlyear}{2014}
\endentry

\entry{djukic_advanced_2015}{inproceedings}{useprefix}
  \name{author}{7}{}{%
    {{}%
     {Djukic}{D.}%
     {Tamara}{T.}%
     {}{}%
     {}{}}%
    {{}%
     {Barcelò}{B.}%
     {Jaume}{J.}%
     {}{}%
     {}{}}%
    {{}%
     {Bullejos}{B.}%
     {Manuel}{M.}%
     {}{}%
     {}{}}%
    {{}%
     {Montero}{M.}%
     {Lidia}{L.}%
     {}{}%
     {}{}}%
    {{}%
     {Cipriani}{C.}%
     {Ernesto}{E.}%
     {}{}%
     {}{}}%
    {{}%
     {Lint}{L.}%
     {Hans}{H.}%
     {van}{v.}%
     {}{}}%
    {{}%
     {Hoogendoorn}{H.}%
     {Serge}{S.}%
     {}{}%
     {}{}}%
  }
  \strng{namehash}{DTBJBM+1}
  \strng{fullhash}{DTBJBMMLCEvLHHS1}
  \field{eventtitle}{Transportation Research Board 94th Annual Meeting}
  \field{shorttitle}{Advanced {{Traffic Data}} for {{Dynamic Origin-Destination
  Demand Estimation}}}
  \field{title}{Advanced Traffic Data for Dynamic Origin-Destination Demand
  Estimation: State of the Art and Benchmark Study}
  \verb{url}
  \verb http://trid.trb.org/view.aspx?id=1338524
  \endverb
  \verb{file}
  \verb Djukic et al_2015_Advanced Traffic Data for Dynamic Origin-Destination
  \verb Demand Estimation.pdf:D\:\\QUT\\Bibliography\\Zotero\\storage\\K5MCRRGR
  \verb \\Djukic et al_2015_Advanced Traffic Data for Dynamic Origin-Destinatio
  \verb n Demand Estimation.pdf:application/pdf;Snapshot:D\:\\QUT\\Bibliography
  \verb \\Zotero\\storage\\PA4I3F5T\\view.html:
  \endverb
  \field{year}{2015}
  \field{urlday}{30}
  \field{urlmonth}{03}
  \field{urlyear}{2015}
\endentry

\entry{perrakis_bayesian_2015}{article}{}
  \name{author}{4}{}{%
    {{}%
     {Perrakis}{P.}%
     {Konstantinos}{K.}%
     {}{}%
     {}{}}%
    {{}%
     {Karlis}{K.}%
     {Dimitris}{D.}%
     {}{}%
     {}{}}%
    {{}%
     {Cools}{C.}%
     {Mario}{M.}%
     {}{}%
     {}{}}%
    {{}%
     {Janssens}{J.}%
     {Davy}{D.}%
     {}{}%
     {}{}}%
  }
  \keyw{Hierarchical Bayesian modelling,Integrated nested Laplace
  approximation,Origin–destination matrix,Overdispersion,Poisson mixtures}
  \strng{namehash}{PKKDCM+1}
  \strng{fullhash}{PKKDCMJD1}
  \field{abstract}{%
  Transportation origin–destination analysis is investigated through the use
  of Poisson mixtures by introducing covariate-based models which incorporate
  different transport modelling phases and also allow for direct probabilistic
  inference on link traffic based on Bayesian predictions. Emphasis is placed
  on the Poisson–inverse Gaussian model as an alternative to the commonly
  used Poisson–gamma and Poisson–log-normal models. We present a first full
  Bayesian formulation and demonstrate that the Poisson–inverse Gaussian
  model is particularly suited for origin–destination analysis because of its
  desirable marginal and hierarchical properties. In addition, the integrated
  nested Laplace approximation is considered as an alternative to Markov chain
  Monte Carlo sampling and the two methodologies are compared under specific
  modelling assumptions. The case-study is based on 2001 Belgian census data
  and focuses on a large, sparsely distributed origin–destination matrix
  containing trip information for 308 Flemish municipalities.%
  }
  \verb{doi}
  \verb 10.1111/rssa.12057
  \endverb
  \field{issn}{1467-985X}
  \field{number}{1}
  \field{pages}{271\bibrangedash 296}
  \field{shortjournal}{J. R. Stat. Soc. A}
  \field{shorttitle}{Bayesian Inference for Transportation Origin–destination
  Matrices}
  \field{title}{Bayesian Inference for Transportation Origin-Destination
  Matrices: The {{Poisson}}-Inverse {{Gaussian}} and Other {{Poisson}}
  Mixtures}
  \field{volume}{178}
  \field{langid}{english}
  \verb{file}
  \verb Snapshot:D\:\\QUT\\Bibliography\\Zotero\\storage\\35RAMCKR\\abstract\;j
  \verb sessionid=524B9ECB8E455C65FCBDDF640B37CC7F.html:;Perrakis et al_2015_Ba
  \verb yesian inference for transportation origin–destination matrices.pdf:D
  \verb \:\\QUT\\Bibliography\\Zotero\\storage\\F59A2SQN\\Perrakis et al_2015_B
  \verb ayesian inference for transportation origin–destination matrices.pdf:
  \verb application/pdf
  \endverb
  \field{journaltitle}{Journal of the Royal Statistical Society: Series A
  (Statistics in Society)}
  \field{year}{2015}
  \field{urlday}{16}
  \field{urlmonth}{01}
  \field{urlyear}{2015}
\endentry

\entry{porta_network_2006}{article}{}
  \name{author}{3}{}{%
    {{}%
     {Porta}{P.}%
     {Sergio}{S.}%
     {}{}%
     {}{}}%
    {{}%
     {Crucitti}{C.}%
     {Paolo}{P.}%
     {}{}%
     {}{}}%
    {{}%
     {Latora}{L.}%
     {Vito}{V.}%
     {}{}%
     {}{}}%
  }
  \keyw{Scale-free networks,Structure of complex networks}
  \strng{namehash}{PSCPLV1}
  \strng{fullhash}{PSCPLV1}
  \field{abstract}{%
  The application of the network approach to the urban case poses several
  questions in terms of how to deal with metric distances, what kind of graph
  representation to use, what kind of measures to investigate, how to deepen
  the correlation between measures of the structure of the network and measures
  of the dynamics on the network, what are the possible contributions from the
  GIS community. In this paper, the author considers six cases of urban street
  networks characterized by different patterns and historical roots. The
  authors propose a representation of the street networks based firstly on a
  primal graph, where intersections are turned into nodes and streets into
  edges. In a second step, a dual graph, where streets are nodes and
  intersections are edges, is constructed by means of a generalization model
  named Intersection Continuity Negotiation, which allows to acknowledge the
  continuity of streets over a plurality of edges. Finally, the authors address
  a comparative study of some structural properties of the dual graphs, seeking
  significant similarities among clusters of cases. A wide set of network
  analysis techniques are implemented over the dual graph: in particular the
  authors show that the absence of any clue of assortativity differentiates
  urban street networks from other non-geographic systems and that most of the
  considered networks have a broad degree distribution typical of scale-free
  networks and exhibit small-world properties as well.%
  }
  \verb{doi}
  \verb 10.1016/j.physa.2005.12.063
  \endverb
  \field{issn}{0378-4371}
  \field{number}{2}
  \field{pages}{853\bibrangedash 866}
  \field{shortjournal}{Physica A: Statistical Mechanics and its Applications}
  \field{shorttitle}{The Network Analysis of Urban Streets}
  \field{title}{The Network Analysis of Urban Streets: A Dual Approach}
  \field{volume}{369}
  \verb{file}
  \verb Porta et al_2006_The network analysis of urban streets.pdf:D\:\\QUT\\Bi
  \verb bliography\\Zotero\\storage\\KGDAFZ53\\Porta et al_2006_The network ana
  \verb lysis of urban streets.pdf:application/pdf;ScienceDirect Snapshot:D\:\\
  \verb QUT\\Bibliography\\Zotero\\storage\\NM87K9F4\\S0378437106001282.html:
  \endverb
  \field{journaltitle}{Physica A: Statistical Mechanics and its Applications}
  \field{day}{15}
  \field{month}{09}
  \field{year}{2006}
  \field{urlday}{17}
  \field{urlmonth}{02}
  \field{urlyear}{2016}
\endentry

\entry{cormen_algorithms_2009}{incollection}{}
  \name{author}{4}{}{%
    {{}%
     {Cormen}{C.}%
     {Thomas~H.}{T.~H.}%
     {}{}%
     {}{}}%
    {{}%
     {Leiserson}{L.}%
     {Charles~E.}{C.~E.}%
     {}{}%
     {}{}}%
    {{}%
     {Rivest}{R.}%
     {Ronald~L.}{R.~L.}%
     {}{}%
     {}{}}%
    {{}%
     {Stein}{S.}%
     {Clifford}{C.}%
     {}{}%
     {}{}}%
  }
  \list{publisher}{1}{%
    {{The MIT Press}}%
  }
  \strng{namehash}{CTHLCERRL+1}
  \strng{fullhash}{CTHLCERRLSC1}
  \field{booktitle}{Introduction to {{Algorithms}}}
  \field{isbn}{978-0-262-03384-8}
  \field{pages}{631\bibrangedash 638}
  \field{title}{The Algorithms of {{Kruskal}} and {{Prim}}}
  \field{year}{2009}
\endentry

\entry{coifman_using_1999}{article}{}
  \name{author}{1}{}{%
    {{}%
     {Coifman}{C.}%
     {Benjamin}{B.}%
     {}{}%
     {}{}}%
  }
  \strng{namehash}{CB1}
  \strng{fullhash}{CB1}
  \field{abstract}{%
  Dual loop speed traps have a distinct advantage over single loop detectors
  because the speed trap detection system is redundant. Each vehicle is
  observed twice under normal operating conditions, once at each loop. The two
  observations are normally used to measure velocity, but as this paper
  demonstrates, the redundancy can also be used to assess the performance of
  the speed trap and identify detector errors. At free-flow velocities, the
  time each detector is occupied by a vehicle (i.e., the ontime) should be
  virtually identical, regardless of the vehicle length. Many hardware errors
  will cause the two on-times to differ. Exploiting this property, a formal
  methodology for testing speed traps off-line has been developed, and ways to
  extend the work to on-line testing are suggested. The work is used to
  evaluate several loop sensor units, revealing problems in two models. A
  second example shows how the work can be used to detect cross talk between
  sensor units.%
  }
  \verb{doi}
  \verb 10.3141/1683-07
  \endverb
  \field{issn}{0361-1981}
  \field{number}{1683}
  \field{pages}{47\bibrangedash 58}
  \field{title}{Using Dual Loop Speed Traps to Identify Detector Errors}
  \verb{file}
  \verb Coifman_1999_Using dual loop speed traps to identify detector errors.pd
  \verb f:D\:\\QUT\\Bibliography\\Zotero\\storage\\IIR99RU9\\Coifman_1999_Using
  \verb  dual loop speed traps to identify detector errors.pdf:application/pdf;
  \verb Snapshot:D\:\\QUT\\Bibliography\\Zotero\\storage\\SG883NEA\\1683-07.htm
  \verb l:
  \endverb
  \field{journaltitle}{Transportation Research Record: Journal of the
  Transportation Research Board}
  \field{year}{1999}
  \field{urlday}{11}
  \field{urlmonth}{02}
  \field{urlyear}{2016}
\endentry

\entry{rubner_earth_2000}{article}{}
  \name{author}{3}{}{%
    {{}%
     {Rubner}{R.}%
     {Yossi}{Y.}%
     {}{}%
     {}{}}%
    {{}%
     {Tomasi}{T.}%
     {Carlo}{C.}%
     {}{}%
     {}{}}%
    {{}%
     {Guibas}{G.}%
     {Leonidas~J.}{L.~J.}%
     {}{}%
     {}{}}%
  }
  \keyw{Artificial Intelligence (incl. Robotics),Automation and
  Robotics,color,Computer Imaging; Graphics and Computer Vision,Earth Mover's
  Distance,image processing,image retrieval,perceptual metrics,texture}
  \strng{namehash}{RYTCGLJ1}
  \strng{fullhash}{RYTCGLJ1}
  \field{abstract}{%
  We investigate the properties of a metric between two distributions, the
  Earth Mover's Distance (EMD), for content-based image retrieval. The EMD is
  based on the minimal cost that must be paid to transform one distribution
  into the other, in a precise sense, and was first proposed for certain vision
  problems by Peleg, Werman, and Rom. For image retrieval, we combine this idea
  with a representation scheme for distributions that is based on vector
  quantization. This combination leads to an image comparison framework that
  often accounts for perceptual similarity better than other previously
  proposed methods. The EMD is based on a solution to the transportation
  problem from linear optimization, for which efficient algorithms are
  available, and also allows naturally for partial matching. It is more robust
  than histogram matching techniques, in that it can operate on variable-length
  representations of the distributions that avoid quantization and other
  binning problems typical of histograms. When used to compare distributions
  with the same overall mass, the EMD is a true metric. In this paper we focus
  on applications to color and texture, and we compare the retrieval
  performance of the EMD with that of other distances.%
  }
  \verb{doi}
  \verb 10.1023/A:1026543900054
  \endverb
  \field{issn}{0920-5691, 1573-1405}
  \field{number}{2}
  \field{pages}{99\bibrangedash 121}
  \field{shortjournal}{International Journal of Computer Vision}
  \field{title}{The {{Earth Mover}}'s Distance as a Metric for Image Retrieval}
  \field{volume}{40}
  \field{langid}{english}
  \verb{file}
  \verb Snapshot:D\:\\QUT\\Bibliography\\Zotero\\storage\\NCHKP6WI\\A1026543900
  \verb 054.html:;Rubner et al_2000_The Earth Mover's Distance as a Metric for
  \verb Image Retrieval.pdf:D\:\\QUT\\Bibliography\\Zotero\\storage\\NJCIWTDN\\
  \verb Rubner et al_2000_The Earth Mover's Distance as a Metric for Image Retr
  \verb ieval.pdf:application/pdf
  \endverb
  \field{journaltitle}{International Journal of Computer Vision}
  \field{month}{11}
  \field{year}{2000}
  \field{urlday}{27}
  \field{urlmonth}{01}
  \field{urlyear}{2016}
\endentry

\entry{combettes_stochastic_2015}{article}{}
  \name{author}{2}{}{%
    {{}%
     {Combettes}{C.}%
     {Patrick~L.}{P.~L.}%
     {}{}%
     {}{}}%
    {{}%
     {Pesquet}{P.}%
     {Jean-Christophe}{J.-C.}%
     {}{}%
     {}{}}%
  }
  \keyw{Mathematics - Optimization and Control,Primary 47H05; Secondary 65K05;
  90C25; 94A08}
  \strng{namehash}{CPLPJC1}
  \strng{fullhash}{CPLPJC1}
  \field{abstract}{%
  We investigate the asymptotic behavior of a stochastic version of the
  forward-backward splitting algorithm for finding a zero of the sum of a
  maximally monotone set-valued operator and a cocoercive operator in Hilbert
  spaces. Our general setting features stochastic approximations of the
  cocoercive operator and stochastic perturbations in the evaluation of the
  resolvents of the set-valued operator. In addition, relaxations and not
  necessarily vanishing proximal parameters are allowed. Weak and strong almost
  sure convergence properties of the iterates is established under mild
  conditions on the underlying stochastic processes. Leveraging these results,
  we also establish the almost sure convergence of the iterates of a stochastic
  variant of a primal-dual proximal splitting method for composite minimization
  problems.%
  }
  \verb{eprint}
  \verb 1507.07095
  \endverb
  \field{title}{Stochastic Approximations and Perturbations in Forward-Backward
  Splitting for Monotone Operators}
  \verb{url}
  \verb http://arxiv.org/abs/1507.07095
  \endverb
  \verb{file}
  \verb arXiv.org Snapshot:D\:\\QUT\\Bibliography\\Zotero\\storage\\3GDP4PPP\\1
  \verb 507.html:;Combettes_Pesquet_2015_Stochastic Approximations and Perturba
  \verb tions in Forward-Backward Splitting for.pdf:D\:\\QUT\\Bibliography\\Zot
  \verb ero\\storage\\NZIKCDXB\\Combettes_Pesquet_2015_Stochastic Approximation
  \verb s and Perturbations in Forward-Backward Splitting for.pdf:application/p
  \verb df
  \endverb
  \field{eprinttype}{arxiv}
  \field{day}{25}
  \field{month}{07}
  \field{year}{2015}
  \field{urlday}{22}
  \field{urlmonth}{02}
  \field{urlyear}{2016}
\endentry

\entry{cremer_new_1987}{article}{}
  \name{author}{2}{}{%
    {{}%
     {Cremer}{C.}%
     {M.}{M.}%
     {}{}%
     {}{}}%
    {{}%
     {Keller}{K.}%
     {H.}{H.}%
     {}{}%
     {}{}}%
  }
  \strng{namehash}{CMKH1}
  \strng{fullhash}{CMKH1}
  \field{abstract}{%
  A new systems dynamics approach for the identification of origin-destination
  (O-D) flows in a traffic system is presented. It is the basic idea of this
  approach that traffic flow through a facility is treated as a dynamic process
  in which the sequences of short-time exit flow counts depend by causal
  relationships upon the time-variable sequences of entrance flow volumes. In
  that way enough information can be obtained from the counts at the entrances
  and the exits to obtain unique and bias-free estimates for the unknown O-D
  flows without further a priori information. Four different methods were
  developed: an ordinary least squares estimator involving cross-correlation
  matrices, a constrained optimization method, a simple recursive estimation
  formula and estimation by Kalman filtering. The methods need only moderate
  computational effort and are particularly useful for tracking time-variable
  O-D patterns for on-line identification and control purposes. An analysis of
  the accuracy of the estimates and a discussion of the convergence properties
  of the methods are given. Finally, a comparison with some conventional static
  estimation procedures is carried out using synthetic as well as real data
  from several intersections. These tests demonstrated that the presented
  dynamic methods are highly superior to conventional techniques and produce
  more accurate results.%
  }
  \verb{doi}
  \verb 10.1016/0191-2615(87)90011-7
  \endverb
  \field{issn}{0191-2615}
  \field{number}{2}
  \field{pages}{117\bibrangedash 132}
  \field{shortjournal}{Transportation Research Part B: Methodological}
  \field{title}{A New Class of Dynamic Methods for the Identification of
  Origin-Destination Flows}
  \field{volume}{21}
  \verb{file}
  \verb Cremer_Keller_1987_A new class of dynamic methods for the identificatio
  \verb n of origin-destination.pdf:D\:\\QUT\\Bibliography\\Zotero\\storage\\PD
  \verb 9QS98P\\Cremer_Keller_1987_A new class of dynamic methods for the ident
  \verb ification of origin-destination.pdf:application/pdf;ScienceDirect Snaps
  \verb hot:D\:\\QUT\\Bibliography\\Zotero\\storage\\SCDWKVDQ\\0191261587900117
  \verb .html:
  \endverb
  \field{journaltitle}{Transportation Research Part B: Methodological}
  \field{month}{04}
  \field{year}{1987}
\endentry

\entry{barcelo_robustness_2013}{book}{}
  \name{author}{5}{}{%
    {{}%
     {Barceló}{B.}%
     {J.}{J.}%
     {}{}%
     {}{}}%
    {{}%
     {Montero}{M.}%
     {L.}{L.}%
     {}{}%
     {}{}}%
    {{}%
     {Bullejos}{B.}%
     {M.}{M.}%
     {}{}%
     {}{}}%
    {{}%
     {Linares}{L.}%
     {M.P.}{M.}%
     {}{}%
     {}{}}%
    {{}%
     {Serch}{S.}%
     {O.}{O.}%
     {}{}%
     {}{}}%
  }
  \strng{namehash}{BJMLBM+1}
  \strng{fullhash}{BJMLBMLMSO1}
  \field{abstract}{%
  Origin-destination (O-D) trip matrices that describe the patterns of traffic
  behavior across a network are the primary data input used in principal
  traffic models and, therefore, a critical requirement in all advanced systems
  supported by dynamic traffic assignment models. However, because O-D matrices
  are not directly observable, the current practice consists of adjusting an
  initial or seed matrix from link flow counts that are provided by an existing
  layout of traffic-counting stations. The availability of new traffic
  measurements provided by information and communication technologies (ICT)
  allows more efficient algorithms, namely for real-time estimation of O-D
  matrices that are based on modified Kalman filtering approaches to exploit
  the new data. The quality of the estimations depends on various factors such
  as the penetration of the ICT devices, the detection layout, and the quality
  of the initial information. The feasibility of real-time applications depends
  on the computational performance of the proposed algorithms for urban
  networks of sensitive size. This paper presents the results of a set of
  computational experiments with a microscopic simulation of the network of
  Barcelona's central business district that explore the sensitivity of the
  Kalman filter estimates in relation to design factor values.%
  }
  \field{pagetotal}{31}
  \field{title}{Robustness and Computational Efficiency of {{Kalman}} Filter
  Estimator of Time-Dependent Origin-Destination Matrices}
  \field{langid}{english}
  \verb{file}
  \verb Barceló et al_2013_Robustness and computational efficiency of kalman f
  \verb ilter estimator of.pdf:D\:\\QUT\\Bibliography\\Zotero\\storage\\ZZNXI3T
  \verb N\\Barceló et al_2013_Robustness and computational efficiency of kalma
  \verb n filter estimator of.pdf:application/pdf
  \endverb
  \field{year}{2013}
\endentry

\entry{lu_kalman_2015}{article}{}
  \name{author}{5}{}{%
    {{}%
     {Lu}{L.}%
     {Zhenbo}{Z.}%
     {}{}%
     {}{}}%
    {{}%
     {Rao}{R.}%
     {Wenming}{W.}%
     {}{}%
     {}{}}%
    {{}%
     {Wu}{W.}%
     {Yao-Jan}{Y.-J.}%
     {}{}%
     {}{}}%
    {{}%
     {Guo}{G.}%
     {Li}{L.}%
     {}{}%
     {}{}}%
    {{}%
     {Xia}{X.}%
     {Jingxin}{J.}%
     {}{}%
     {}{}}%
  }
  \keyw{dynamic OD estimation,dynamic traffic assignment,Kalman filter,Traffic
  sensors,traffic simulation,urban network}
  \strng{namehash}{LZRWWYJ+1}
  \strng{fullhash}{LZRWWYJGLXJ1}
  \field{abstract}{%
  Considerable efforts have been devoted to the development of dynamic
  origin-destination (OD) estimation models, which are a key step to realizing
  self-adaptive traffic control systems for urban traffic management. However,
  most of the models proposed to date estimate OD flows based on a single
  traffic data source, and their performance is limited by the coverage and
  accuracy of traffic sensors. The inherent difficulty in estimating the
  dynamic traffic assignment matrix means that dynamic OD estimation remains a
  challenge for real-life applications. This paper proposes the use of a Kalman
  filter for dynamic OD estimation using multi-source sensor data. The dynamic
  characteristic of changing OD flow over time is analyzed, and the problem of
  dynamic OD estimation is converted to a problem of estimating OD structural
  deviation. The resulting dynamic relationship between traffic volume and OD
  structural deviation is then used to establish the Kalman filter model. An
  improved traffic assignment approach is developed and embedded into the
  measurement equation of the Kalman filter model to enable dynamic updating of
  the traffic assignment matrix. A dual self-adaptive mechanism based on the
  Kalman filter is used to calibrate the model. The proposed method was
  implemented on a real-life traffic network in the downtown area of Kunshan
  City, China. The results show that the proposed method is more accurate than,
  and outperforms, the traditional link-volume-based and turning-movement-based
  methods. Copyright © 2014 John Wiley \& Sons, Ltd.%
  }
  \verb{doi}
  \verb 10.1002/atr.1292
  \endverb
  \field{issn}{2042-3195}
  \field{number}{2}
  \field{pages}{210\bibrangedash 227}
  \field{shortjournal}{J. Adv. Transp.}
  \field{title}{A {{Kalman}} Filter Approach to Dynamic od Flow Estimation for
  Urban Road Networks Using Multi-Sensor Data}
  \field{volume}{49}
  \field{langid}{english}
  \verb{file}
  \verb Snapshot:D\:\\QUT\\Bibliography\\Zotero\\storage\\7JPSDHR7\\abstract\;j
  \verb sessionid=38FF4EB1CFEB53D2C267C07C93DE907F.html:;Lu et al_2015_A Kalman
  \verb  filter approach to dynamic OD flow estimation for urban road networks.
  \verb pdf:D\:\\QUT\\Bibliography\\Zotero\\storage\\T6IK92B5\\Lu et al_2015_A
  \verb Kalman filter approach to dynamic OD flow estimation for urban road net
  \verb works.pdf:application/pdf
  \endverb
  \field{journaltitle}{Journal of Advanced Transportation}
  \field{day}{01}
  \field{month}{03}
  \field{year}{2015}
  \field{urlday}{29}
  \field{urlmonth}{01}
  \field{urlyear}{2016}
\endentry

\lossort
\endlossort

  \blx@bblend
  \endgroup
  \csnumgdef{blx@labelnumber@\the\c@refsection}{0}}
\makeatother

\usepackage{amssymb}
\usepackage{color}

\usepackage{colortbl}
\definecolor{lgray}{gray}{0.75}
\definecolor{dgray}{gray}{0.55}

\usepackage{authblk}

\DeclareGraphicsExtensions{.jpg, .pdf, .png}
\graphicspath{{Figures/}}

\newcommand{\udvar}[1]{\underline{{#1}}}
\newcommand{\ddvar}[1]{\underline{\underline{{#1}}}}
\newcommand{\tdvar}[1]{\underline{\underline{\underline{{#1}}}}}
\newcommand{\qdvar}[1]{\underline{\underline{\underline{\underline{#1}}}}}

\newcommand{\RR}{\ensuremath{\mathbb R}}

\newcommand{\esp}{\ensuremath{\mathbb E }}
\newcommand{\prox}{\ensuremath{\mathrm{prox}}}
\newcommand{\Argmin}{\ensuremath{\mathrm{Argmin}\,}}

\newcommand{\sign}{\ensuremath{\mathrm{sign}}}

\def\eg{e.g.,\ }
\def\cf{\textit{cf.\ }}
\def\ie{\textit{i.e.,\ }}

\begin{document}

\title{A Primal-Dual Algorithm for Link Dependent Origin Destination Matrix Estimation\thanks{Preliminary versions of this work were presented in \cite{michau_estimating_2015-1,michau_estimating_2015}. Work supported by ANR-14-CE27-0001 GRAPHSIP grant and ANR-12-SOIN-0001-02 V\'el'Innov grant.}}

\author[$\times$,$\dag$]{Gabriel Michau}
\author[$\times$]{Nelly Pustelnik}
\author[$\times$]{Pierre Borgnat}
\author[$\times$]{Patrice Abry}
\author[$\dag$]{Alfredo Nantes}
\author[$\dag$]{Ashish Bhaskar}
\author[$\dag$]{Edward Chung}
\affil[$\times$]{Univ Lyon, Ens de Lyon, Univ Lyon 1, CNRS, Laboratoire de Physique, F-69342 Lyon, France}
\affil[$\dag$]{Queensland University of Technology, Smart Transport Research Centre, Brisbane, Australia}

\maketitle

\begin{abstract}
Origin-Destination Matrix (ODM) estimation is a classical problem in transport engineering 
aiming to recover flows from every Origin to every Destination from measured traffic counts and a priori model information.
In addition to traffic counts, the present contribution takes advantage of probe trajectories, whose capture is made possible by new measurement technologies. It extends the concept of ODM to that of Link dependent ODM (LODM), keeping the information about the flow distribution on links and containing inherently the ODM assignment.
Further, an original formulation of LODM estimation, from traffic counts and probe trajectories is presented  
as an optimisation problem, where the functional to be minimized consists of five convex functions, each modelling a constraint or property of the transport problem: 
consistency with traffic counts, 
consistency with sampled probe trajectories, 
consistency with traffic conservation (Kirchhoff's law), 
similarity of flows having close origins and destinations, 
positivity of traffic flows. 
A primal-dual algorithm is devised to minimize the designed functional, as the corresponding objective functions are not necessarily differentiable. 
A case study, on a simulated network and traffic, validates the feasibility of the procedure and details its benefits for the estimation of an LODM matching real-network constraints and observations.
\end{abstract}

\section{Introduction}

The estimation of traffic flows over networks is a keystone for understanding their usage and behaviour in specific situations, e.g., network has a limited capacity or traffic may significantly vary with time or with particular events. 
Estimating traffic flows is thus needed for the network efficiency analysis, for traffic prediction, and traffic optimisation.
Origin-Destination matrices (ODM) estimation is one of the classical problem in transport engineering \cite{willumsen_estimation_1978} but also in the study of Internet traffic \cite{coates_internet_2002,girard_performance_2006,mardani_estimating_2015}.
ODM are double entry tables indexed by network zones or major nodes, whose elements contain the demand of traffic from origins indexed by rows, to destinations, indexed by columns. With respect to the transport field, ODM can be recovered from traveller interviews directly. This is however a long, difficult and costly process. Thus, since the 70's and as a consequence of the generalisation, in occidental cities, of the access to link counts (\eg by magnetic/inductive loops), many researches have sought to estimate the ODM with traffic counts as their primary source of data.
	
\noindent \textbf{Estimating ODM from link counts.} \quad 
Formally, a road network is represented by a graph $\mathcal{G} = (V,L)$, with $\ddvar{T}$ the corresponding ODM, of size $ \vert V\vert\times \vert V\vert$. 
Magnetic loops, on links $l\in L$, produce $ \vert L\vert$ measures represented by vector $\udvar{q}$. 
Thus, ODM estimation problem amounts to solving the following inverse problem:
\begin{equation}
\label{eqn:intro:invPb}
 \udvar{q} = F(\ddvar{T}) +  \udvar{\epsilon}
\end{equation}
where the assignment function $F$ relates OD flows to network link, for comparisons against traffic counts $q$, and where $\udvar{\epsilon}$ models 
the measurement error.
The two main difficulties in solving Problem \eqref{eqn:intro:invPb} stem, first, from its being ill-conditioned:  the size of the quantity to be estimated $\ddvar{T}$ is larger than that of the available measures $\udvar{q}$ and second, from $F$ being unknown and thus often modelled.

To solve Problem \eqref{eqn:intro:invPb}, a common approach is to rely on the so-called four-steps model, that consists of Trip Generation, Trip Distribution, Modal Split\footnote{The Modal Split's interest lies when one consider several modes of transport. Here however, the focus is on car trips only and this step is ignored.}, and Trip Assignment. 
The first two steps permit to design $\ddvar{T}$ while the Assignment step, amounts to specifying $F$.  
Interested readers can refer to \cite{ortuzar_modelling_2011} for more information on this framework.	
The use of this model requires a fine parameters tuning for the three steps. 
Hence, a large literature can be found, with numerous variations for both trip distribution and assignment. 
A detailed review is beyond the scope of the present contribution, and interested readers are referred to e.g., \cite{wilson_entropy_1970, willumsen_estimation_1978, parry_estimation_2012,kim_estimation_2010,castillo_bayesian_2013,mellegard_origin/destination-estimation_2011, iqbal_development_2014,alexander_origin-destination_2015} and references therein.\\

\noindent \textbf{Goals, contributions and outline.} \quad 
Despite the fact that  $\ddvar{T}$ is of size $\vert V\vert \times \vert V\vert$, solving \eqref{eqn:intro:invPb} is in fact an inverse problem of size $\vert V\vert \times \vert V\vert \times \vert L\vert$, 
because of the required assignment step that actually involves the number of links in the network.
The goal of the present contribution is to directly account for the real dimensionality of the problem by proposing a new and original description tool for traffic, that directly includes assignment: the Link dependent Origin Destination Matrix (LODM). 
LODM represents the OD flows already assigned to each link of the network, thus incorporating the assignment, or equivalently making its independent specification unnecessary. 
We also propose to estimate LODM as an inverse problem of dimension $\vert V\vert \times \vert V\vert \times \vert L\vert$. We rely on traffic counts $q$ and, in addition, on a new set of data: partial (or sampled) knowledge of trajectories, whose collection is now made possible by new technologies such as GPS \cite{herrera_evaluation_2010}, Bluetooth \cite{michau_bluetooth_2016, hainen_estimating_2011, feng_vehicle_2015}, Floating car data \cite{gomez_evaluation_2015}.
Section \ref{sec:RoadNetwork} formalises the transport problem, from its engineering perspective.
Section \ref{sec:VA} details five significant properties imposed either by the network or for consistency with the observed data and turns them into five components of an objective function that formalises LODM estimation from traffic counts and sampled trajectories.
Because these five functions are convex but non necessarily differentiable, a proximal primal-dual algorithm is devised to minimize the corresponding optimisation problem results.
To finish with, the feasibility of the proposed approach and the assessment of its estimation performance are investigated in Section \ref{sec:SCS}, on a case study consisting of network and traffic simulations, designed to match closely various realistic network and traffic in large western metropolitan cities.\\
	
\noindent \textbf{Notations.} \quad 
	The following notations are used throughout this article: 
	$\udvar{X}$, $\ddvar{X}$ and $\tdvar{X}$ refer to vectors, matrices and tensors, respectively.
The Hadamard product (element-wise product) of $\tdvar{Y}$ and $\tdvar{X}$ is denoted $\tdvar{Y}\circ \tdvar{X}$.
Subscript indices are used for dimensions over the nodes of the graph and the index $i$ is used to label origins, $j$ to label destinations, $k, m, n$ and $p$ to label nodes in general. Superscript indices are used for dimensions over the links and the indexes $l$ and $e$ are favoured.

The symbol $\bullet$ is used to denote the dimension that does not contribute to a sum: \eg the sum over first and third dimensions, indexed respectively with $i$ and $l$, is written $\sum_{i\bullet l} \tdvar{X}$.

We denote by $\Vert \cdot \Vert_1$ the element-wise first norm for matrices: \eg $\Vert \ddvar{X} \Vert_1 = \sum_{ij} \vert X_{ij} \vert$.

\section{Road Network and Link-Dependent ODM}
\label{sec:RoadNetwork}

	\subsection{The problem}
\begin{figure}[t]
	\centering
	\subfloat[]{\includegraphics[width=5cm]{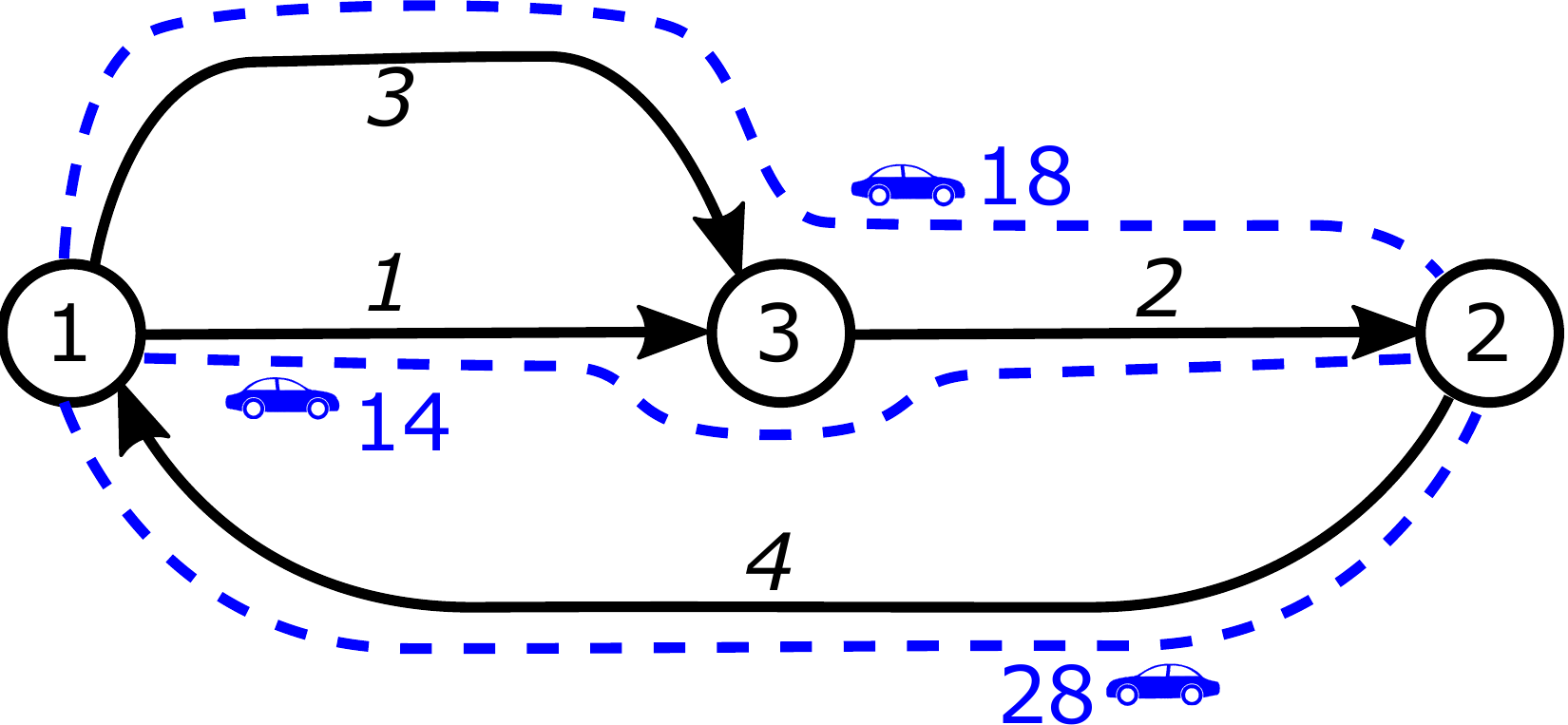}\label{sfig:Network_Example}}\\
  	\subfloat[]{\renewcommand{\arraystretch}{0.5} \renewcommand{\tabcolsep}{1.5pt}
  	\shortstack{
  	$ \ddvar{I} = \left[ 
  	\begin{tabular}{c c c c}
  	0&0&0&1\\0&1&0&0\\1&0&1&0
  	\end{tabular} \right] $ \hfil
  	$ \ddvar{E} = \left[ 
  	\begin{tabular}{c c c c}
  	1&0&1&0\\0&0&0&1\\0&1&0&0
  	\end{tabular} \right] $ \hfil
  	$ \udvar{q} = \left[ 
  	\begin{tabular}{c}
  	14\\32\\18\\28
  	\end{tabular}  \right] $ \hfil
  	  	$ \ddvar{T} = \left[ 
  	\begin{tabular}{c c c}
  	0&32&0\\28&0&0\\0&0&0
  	\end{tabular} \right] $ 
	\\ \raisebox{1cm}{$\tdvar{Q} =\ \ $}\includegraphics[width=4cm]{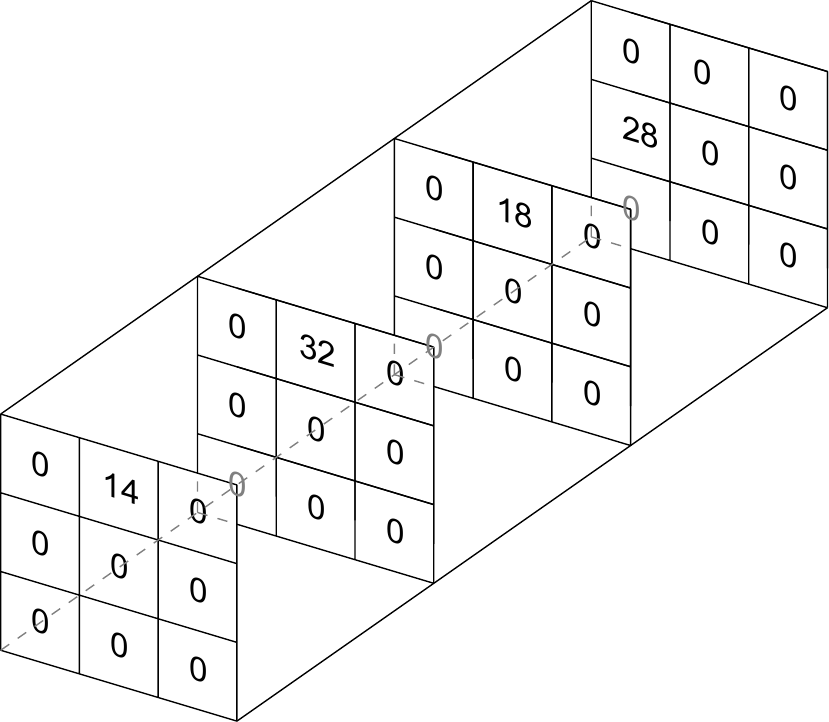}	}\label{sfig:LSODM}}
	\caption{Example of a simple network (a) with the associated tools describing the topology and the traffic (b).}
	\label{fig:Network_Example}
\end{figure}

	The network is described as a graph $\mathcal{G}=(V,L)$ where the finite set of nodes $V$ models intersections of the road network. 
	Each node also defines a possible origin or destination. 
$L$ is the set of directed edges, each corresponding to a direct itinerary (or road) linking two nodes (\ie not going through another node in $V$). 
The number of road users is denoted $N$. 
A schematic (small) such graph is illustrated in Fig. \ref{sfig:Network_Example}.

On such a graph, LODM consists of a tensor of size $ \vert V\vert \times \vert V\vert \times \vert L\vert$, labelled $\tdvar{Q} = (Q_{ij}^l)_{{(i,j)\in V^2, l \in L}}.$ 
As illustrated in Fig. \ref{fig:Network_Example},  each trajectory adds a count of 1 in $Q_{ij}^l$ if the link $l$ is on the origin-destination path $(i,j)$. 
Therefore, $\tdvar{Q}$ consists, for each link $l\in L$, in an OD matrix of size $\vert V\vert \times \vert V\vert$.

To perform the estimation of $\tdvar{Q}$, we use information stemming from probe trajectories as well as traffic counts on each link. 
The set of trajectories can be measured from various sources (GPS, Bluetooth, ...) and the actual technology matters little in the procedure. We propose here, though, to refer to the Bluetooth technology, which is of great interest as it currently provides trajectory datasets with the highest penetration rate, compared to other technologies; the penetration rate is the fraction of vehicles equipped with the chosen technology and from which information needed to reconstruct trajectory can be collected \cite{michau_bluetooth_2016}.
Trajectory information is stored into a tensor labelled $\tdvar{B}$, of size $|V|\times|V|\times |L|$. 
This tensor can be read as a sampled version of $\tdvar{Q}$, from only a fraction of the total traffic. 
Traffic counts consists of the total volume of traffic on each link $l\in L$, labelled $\udvar{q}$, of size $\vert L\vert$, irrespective of OD pairs. 
Traffic counts can be, for instance, measured by magnetic loops.

A variational approach will now be devised to estimate $\tdvar{Q}^\star$, the real LODM, by means of non-smooth convex optimisation from $\tdvar{B}$ and $\udvar{q}$. 
The involved criterion represents on the one hand the relationship between the tensor $\tdvar{Q}$ and the measures ($\tdvar{B}$, $\udvar{q}$) and, on the other hand, properties of the road network and traffic constraints (e.g., car conservation at intersections).

	\subsection{Structure of the graph and of the traffic}
The structure of the graph is given by the \textit{incidence} and \textit{excidence} matrices denoted respectively $\ddvar{I}$ and $\ddvar{E}$ of size $\vert V\vert \times \vert L \vert$. 
These matrices describe the relations between the nodes and the edges, such that, for every $(k,l)\in V \times L$,
\begin{equation}
	\begin{array}{ll}
	 I^l_k\ & = \left\lbrace 
		\begin{array}{ll}
		1 & \mbox{if the link $l$ is arriving to the node $k$}, \\
		0  & \mbox{otherwise},\\
		\end{array}
		\right. \\
		\\
	 E^l_k\ & = \left\lbrace 
		\begin{array}{ll}
		1 & \mbox{if the link $l$ is starting from the node $k$}, \\
		0  & \mbox{otherwise}.\\
		\end{array}
		\right. \\
	\end{array}
\end{equation}

Note that in graph theory, it is customary to name the difference $(\ddvar{I} - \ddvar{E})$ as \textit{Incidence Matrix}; 
however we need both matrices separately in this work.

Let us also define the tensors $\tdvar{I_1}$ and $\tdvar{I_2}$ (resp. $\tdvar{E_1}$ and $\tdvar{E_2}$) corresponding to the replication of $\ddvar{I}$ (resp. $\ddvar{E}$) such that,
\begin{equation}
	\begin{array}{rl}
		(\forall m \in V) \;  & (I_1)^l_{km} =  \left\lbrace
		\begin{array}{l l}
		1 & \mbox{if link $l$ is arriving to node $k$,}\\
		0 & \mbox{otherwise,}\\
		\end{array}
		\right.\\
		(\forall k \in V) \;  & (I_2)^l_{km} =  \left\lbrace
		\begin{array}{l l}
		1 & \mbox{if link $l$ is arriving to node $m$,}\\
		0 & \mbox{otherwise.} \\
		\end{array}
		\right.\\
	\end{array}
\end{equation}

Using these notations, we relate the LODM $\tdvar{Q}$ to the classical OD matrix $\ddvar{T}$ of size $\vert V\vert\times \vert V\vert$ where each element $T_{ij}$ contains the traffic flow originating from the node $i$ and having $j$ for destination as follows:

\begin{equation}
	\label{eqn:ODM}
	\ddvar{T} = \sum_{\bullet\bullet l} \tdvar{E_1}\circ \tdvar{Q} = \sum_{\bullet \bullet l} \tdvar{I_2}\circ \tdvar{Q}.
\end{equation}

We denote by $\udvar{O}$ (resp. $\udvar{D}$) the \textit{origin} (resp. \textit{destination}) vector, of size $\vert V\vert$ as the sum of $\ddvar{T}$ over the second (resp. first) dimension. It represents the flows originating (or having for destination) each node of the graph. Formally,

\begin{equation}
	\label{eqn:ODVec}
	\left\lbrace
		\begin{array}{l}
		\udvar{D} = (\sum_{i\bullet} \ddvar{T})^\top, \\
		\udvar{O} = \sum_{\bullet j} \ddvar{T}. 
		\end{array}
	\right.
\end{equation}

	\subsection{Model, Measures and Estimates}
For this problem, we consider an urban road network composed of major roads, ignoring residential and service streets, seldom equipped with traffic sensors.

The set of users with their trajectories on those major roads, are represented through the tensor $\tdvar{Q}$, as described above and that we wish to estimate.

For this estimation, we first assume that every road is equipped with a magnetic loop, counting the number of cars using it. It implies therefore that every element in $\udvar{{q}}$ is known. This assumption is realistic in our case considering major roads only. 
The magnetic loops are usually subject to counting errors and it is modelled here by a noise $\udvar{\varepsilon}$. Hence the measured quantity $\udvar{q}$ reads:
\begin{equation}
	\label{eqn:TrafficAggr}
	\udvar{{q}} = \udvar{q}^\star + \udvar{\varepsilon}
\end{equation}
where $\udvar{q}^\star$ is the true traffic volumes.

Second, we also assume that the proportion of Bluetooth equipped vehicle can vary with each possible pair of OD but that it doesn't change along a trajectory (the tracking devices are not turned off and on while the car is running). This assumption allows us to define an OD-dependant penetration rate $\ddvar{\eta_o}$ of size $\vert V \vert \times\vert V \vert$.  Experiments in Brisbane have shown that the average Bluetooth penetration rate is around 25\% \cite{laharotte_spatiotemporal_2015}. Moreover, $\tdvar{B}$ appears as a noisy version of $\tdvar{Q}^\star$ for which the noise level depends on the penetration rate. The relation between the tensors $\tdvar{B}$ to $\tdvar{Q}^\star$ can thus be modelled by a Poisson law, typically involved in counting processes. This leads to a model
\begin{equation}
	\label{eqn:Poisson}
	(\forall i,j,l \in V \times V \times L) \quad  B_{ij}^l = \mathcal{P} ((\eta_o)_{ij} Q_{ij}^{\star l}).
\end{equation}

\section{Variational Approach}
\label{sec:VA}
Instead of using the traditional four-steps model resolution, iterating over a process involving \textit{a priori} information, modelling of the traffic, estimating the variables of interest, comparing to the observed measures and tuning the models, we propose here the use of a variational approach. Both our knowledge of the network and of the traffic states are included within an objective function that combines together five terms to be jointly minimised.

	\subsection{Objective Function}

The terms of the objective functions can be classified in three types: The first type, composed of functions \ref{sssec:TrafficCount}, \ref{sssec:BluetoothPoisson} and \ref{sssec:DefinitionDomain}, is aiming for consistency between the measures and the estimate. The second type, with function \ref{sssec:Kirchhoff}, stems from the topology of the network. The third and last type, with function \ref{sssec:TV}, comes from an additional assumption based on our knowledge of transport networks.

		\subsubsection{Traffic Count Data Fidelity $f_{TC}$}
		\label{sssec:TrafficCount}
Ensuring the consistency with traffic counts would require that Eq. \eqref{eqn:TrafficAggr} is satisfied. Noting that:
\begin{equation}
	\label{eqn:TrafficAggr-2}
	{\udvar{q}^\star} = \sum_{ij\bullet} \tdvar{Q}^\star
\end{equation}
and assuming a random unbiased Gaussian noise $\udvar{\varepsilon}$ for the magnetic loops, as in Eq. \eqref{eqn:TrafficAggr}, the constraint of Eq \eqref{eqn:TrafficAggr-2} can be released and leads to the following function:

\begin{equation}
	\label{eqn:Obj:TrafficCount}
	f_{TC}(\tdvar{Q}) = \Vert \udvar{{q}}  -  \sum_{ij\bullet} \tdvar{Q}\Vert^2.
\end{equation}

		\subsubsection{Poisson Bluetooth Sampling Data Fidelity $f_P$}
		\label{sssec:BluetoothPoisson}
Second, the consistency with Bluetooth measures, as modelled in Equation \eqref{eqn:Poisson} requires the knowledge of the OD-dependent penetration rate $\ddvar{\eta_o}$. 
This information, of size  $\vert V\vert \times \vert V\vert$, is not directly available from $\udvar{q}$ and $\tdvar{B}$, therefore we introduce an approximation of this penetration rate of size $\vert L\vert$, noted $\udvar{\eta}$ and calculated as:

\begin{equation}
	\label{eqn:Etal}
	\udvar{\eta} = \frac{\udvar{{q}}}{\sum_{i,j,\bullet} \tdvar{B}}.
\end{equation}

The resulting data fidelity term, denoted $f_P$, models the minus log-likelihood associated with the Poisson model \cite{combettes_douglas-rachford_2007}:
\begin{equation}
	\label{eqn:Obj:Poisson}
	f_P(\tdvar{Q}) = \sum_{ijl}\psi\left(B_{ij}^l,\eta^l Q_{ij}^l \right)
\end{equation}
where, for every $(u,v)\in \RR^2$,
\begin{equation}
	\label{eqn:psi}
	\psi(u, v) =  
	\begin{cases}
	-u \log  v +  v  & \mbox{if}\;\; v>0 \;\;\mbox{and}\;\; u>0,\\
	  v & \mbox{if}\;\; v\geq 0 \;\;\mbox{and}\;\; u=0,\\
	 +\infty & \mbox{otherwise}.\\
	\end{cases}
\end{equation}

		\subsubsection{Definition Domain Constraint $f_C$}
		\label{sssec:DefinitionDomain}
Third, another term ensuring data consistency models that the total flow should be greater than the flow of Bluetooth enabled vehicles. It consists thus in imposing that  $\tdvar{Q}$ belongs to the following convex set $C$:
\begin{equation}
	\label{eqn:ConvexSet}
	C = \big\{\tdvar{Q}  = \big({Q}_{ij}^l\big)_{(ijl)\in V\times V\times L} \in \RR^{\vert V\vert \times \vert V\vert \times \vert L\vert}\;\vert\;  {Q}_{ij}^l\geq {B}_{ij}^l\big\}.
\end{equation}
The corresponding convex function is the indicator function $\iota_C$:
\begin{equation}
	\label{eqn:Obj:domain}
	f_C(\tdvar{Q}) = \iota_C(\tdvar{Q}) = 
	\begin{cases}
	0 &\mbox{\ if\ } \;\;\tdvar{Q}\in C, \\
	+\infty &\mbox{\ otherwise.}
	\end{cases}
\end{equation}

		\subsubsection{Kirchhoff's Law $f_K$}
		\label{sssec:Kirchhoff}
This property is the classical law for flows on network, the Kirchhoff's law, describing the conservation of cars at intersections. It takes into account the network topology. It requires that, for each OD pairs and at every node, the number of cars is conserved when properly accounting for origins and destinations. For every origin $i \in V$, destination $j \in V$ and node $k\in V$ of the network, this yields to, 
\begin{equation}
	\label{eqn:Kirchhoff2}
	\sum_l E^l_k Q^l_{ij}  -  \underbrace{\delta_{ik} T_{ij}}_{\substack{\textrm{origin}\\ \textrm{(source)}}} = \sum_l  I^l_k Q_{ij}^l- \underbrace{\delta_{jk} T_{ij}}_{\substack{\textrm{destination}\\ \textrm{(sink)}}}.
\end{equation}
This constraint can then be summarized as
\begin{equation}
	\label{eqn:Kirchhoff-Q2}
	\big( \forall (i,j,k) \in V^3\big) \quad \sum_{l} A^l_{ijk} Q^l_{ij} = 0
\end{equation}
where the $\vert V\vert   \times \vert V\vert  \times \vert V\vert\times \vert L\vert$  tensor $\qdvar{A}$ is defined as
\begin{equation}
	\label{eqn:Kirchhoff-Op}
	\big(\forall (i,j,k,l) \in V^3\times L\big) \quad 
	A^l_{kij} = \left(E^l_k-I^l_k\right)-\left(\delta_{ik}-\delta_{jk}\right)E^l_i.
\end{equation}
It results in a convex function to be minimised:
\begin{equation}
	\label{eqn:Obj:Kirchhoff}
	f_K(\tdvar{Q}) = \sum_{ijk}\Big(\sum_{l} A^l_{ijk} Q^l_{ij} \Big)^2.
\end{equation}

Compared to our previous works \cite{michau_estimating_2015-1,michau_estimating_2015}, here the Kirchhoff's law is applied per OD pairs and not simply at a global scale. Indeed, the Kirchhoff's law needs also to be satisfied at each node, independently of the origin and destination of the cars. Satisfying this global Kirchhoff's law of \cite{michau_estimating_2015-1,michau_estimating_2015} is a consequence of the one used here. 

		\subsubsection{Total Variation $f_{TV}$}
		\label{sssec:TV}
Finally, from a transport perspective it seems realistic to assume that for two paths having close origins (resp. destinations) and same destination (origin), the trajectories in the network should be correlated (\eg use of similar roads). Such property can be written as
\begin{align}
(\forall i\sim i')(\forall j \in V)(\forall l\in L) \qquad  Q_{ij}^l \sim Q_{i'j}^l \\
(\forall j\sim j')(\forall i \in V)(\forall l\in L) \qquad  Q_{ij}^l \sim Q_{ij'}^l
\end{align}
In order to be used in a variational approach these relationships can be gathered within the convex function $f_{TV}$ defined as the total variation:
\begin{equation}
\label{eqn:TVobj}
f_{TV}(\tdvar{Q} ) =\sum_{i\sim \mathcal{N}_{i'}}\sum_{j,l}   \omega_{ii'} \vert Q_{ij}^l - Q_{i'j}^l\vert + \sum_{j\sim \mathcal{N}_{j'}}\sum_{i,l}   \omega_{jj'} \vert Q_{ij}^l - Q_{ij'}^l\vert 
\end{equation}
where $\mathcal{N}_{i'}$ models the neighbourhood of $i'$ and where $\omega_{ii'}$ are positive weights on edges detailed later.
The use of the $\ell_1$-norm is justified for its \textit{edge preservation} properties. Indeed, it has been shown in \cite{chambolle_introduction_2010, couprie_dual_2013} that the $\ell_1$-norm is adapted for cases where one seeks for spatial correlations while allowing some irregularities, \eg \textit{edges}, in image analysis. From a traffic perspective, we want to encourage users from similar origin (resp. destination) and with same destination (resp. origin) to use similar routes, but also want to allow some irregularities, \eg for nodes in between two major roads where both could be a possible choice. Those nodes can be interpreted as \textit{edges} in image analysis.
 
Equation \eqref{eqn:TVobj} can further be simplified using a weighted effective incidence matrix, denoted $\ddvar{J}$, defined as
\begin{equation}
\big(\forall (k,l)\in V\times L\big)\qquad J_{k}^l = W^l(I_{k}^l - E_{k}^l) 
\end{equation}
and thus having a size $\vert V \vert \times \vert L \vert $, where each element $W^l$ denotes the weight for the link $l$. For this work we choose the following vector $\udvar{W}$ of size $\vert L \vert $:
\begin{equation}
\big(\forall l\in L\big)\qquad W^l = e^{-\frac{d_{l}}{d_0}}
\end{equation}
where $d_l$ models the length of the link $l$ and $d_0$ is the average distance of the nodes.
For the simulated network:
\begin{equation}
d_0 = \sqrt{\frac{\mbox{GridWidth}\cdot\mbox{GridHeight}}{\vert V\vert}}.
\end{equation}
Note that $ W^l = \omega_{km}$ if $l$ is a link between $k$ and $m$.

Eq. \eqref{eqn:TVobj} can then be rewritten as 
\begin{equation}
\label{eqn:TVobj2}
f_{TV}(\tdvar{Q} ) = \sum_{l}   \Vert \ddvar{J}^\top   \ddvar{Q}^l \Vert_1 + \sum_{l}   \Vert \ddvar{J}^\top   (\ddvar{Q}^l)^\top \Vert_1.
\end{equation}
where $\ddvar{Q}^l$ models the $l$-th extracted matrix from $\tdvar{Q}$. Its dimension is thus $\vert V \vert \times \vert V \vert $.

	\subsection{Algorithm}
	\label{ssec:Algo}

To sum up, the objective is to find an estimate of $\tdvar{Q}^\star$ satisfying
\begin{multline}
\label{eqn:Obj}
{\tdvar{\widehat Q}}  \in \underset{\tdvar{Q}}{\Argmin} \big\{\gamma_{TC} f_{TC} + \gamma_{P} f_P + \gamma_{C} f_C \\+ \gamma_{K}f_{K} + \gamma_{TV} f_{TV}\big\}
\end{multline}
where $\gamma_\cdot$ are positive weights applied to the objectives and model their relative importance within the global objective.

All the five functions involved in Eq. \eqref{eqn:Obj} follow the usual assumptions required when dealing with convex optimisation tools: they are convex, lower-semicontinuous (l.s.c.) and proper. Moreover, both the functions $f_{TC}$ and $f_K$ are differentiable and their gradients are given below:
\begin{equation}
\nabla f_{TC}(\tdvar{Q}) = \left(\left( -2 \left( q^l - \sum_{k,m} Q_{km}^l \right) \right)_{ij}^l \right)_{(ijl)\in V\times V\times L}
\end{equation}
and
\begin{equation}
\nabla f_K(\tdvar{Q}) = \left(\left( 2 \sum_k A^l_{jik} \  \sum_e A^e_{ijk} \  Q^e_{ij}) \right)^l_{ij}\right)_{(ijl)\in V\times V\times L}.
\end{equation}
Their Lipschitz constants are denoted $\beta_{TC}$ and $\beta_K$ respectively \cite{pustelnik_wavelet-based_2016}. The other three functions however are not differentiable and $f_{TV}$ involves a linear transformation $H$ such as:
\begin{equation}
	\label{eqn:Obj:FTV}
	f_{TV}(\tdvar{Q}) = \Vert H(\tdvar{Q}) \Vert_1
\end{equation}
where $H$ satisfies:
\begin{equation}
\begin{aligned}
	H\colon & \RR^{\vert V \vert \times \vert V \vert \times \vert L \vert} & \to &\ \RR^{\vert L \vert \times \vert V \vert \times \vert L \vert}\times  \RR^{\vert L \vert \times \vert V \vert \times \vert L \vert}\\
	&\tdvar{Q} & \mapsto & \ \Big(\big( \ddvar{J}^\top   \ddvar{Q}^l\big)_{l\in L}, \big( \ddvar{J}^\top   (\ddvar{Q}^l)^\top \big)_{l\in L}\Big)
\end{aligned}
\end{equation}
and whose adjoint is
\begin{align}
\label{eqn:Hstar}
H^*\colon & (\tdvar{R},\tdvar{S}) & \mapsto & \ \Big( \ddvar{J} \;\ddvar{R}^l\Big)_{l\in L} + \Big( (\ddvar{J} \;\ddvar{S}^l)^\top\Big)_{l\in L}.
\end{align}
In the following, we denote $\chi_{H}$ the norm of this operator. For  further details about the way to compute this norm, the reader can refer to \cite{pustelnik_wavelet-based_2016}.

\begin{algorithm*}[t]
\textbf{Choose:} $\gamma_{TC} \geq 0, \quad \gamma_K \geq 0, \quad \gamma_{TV} \geq 0, \quad \gamma_P\in [0,1], \quad \gamma_D\in [0,1]$\\
\textbf{Compute:} $\chi_H, \quad \beta = \gamma_{TC}\beta_{TC} +\gamma_{K}\beta_{K} $, \quad if $\gamma_{TV} = 0$, $ \quad \tau = \frac{1.99}{\beta}$, else choose $(\tau, \sigma)$ such as $\tau=\frac{0.9}{\frac{\beta}{2}+\sigma \chi_H} \in [\frac{2}{3}\sigma, \frac{3}{2}\sigma]$\\
\textbf{Set :} $\tdvar{Q}^0  = \tdvar{0}$, $ \quad (\tdvar{R}^{0},\tdvar{S}^{0})   = (\tdvar{0}, \tdvar{0})$

\textbf{For} $k = 0$,... :\\
\hspace*{0.4cm} 1. $\widetilde{\tdvar{Q}}^{k+1} =\tdvar{Q}^k - \tau \big(\gamma_{TC} \nabla f_{TC}(\tdvar{Q}^k) + \gamma_{K} \nabla f_{K}(\tdvar{Q}^k) \big) - \tau H^* (\tdvar{R}^k,\tdvar{S}^k ) $\\
\hspace*{0.4cm} 2. $\tdvar{Q}^{k+1} = \prox_{\gamma_C f_C} \Bigg( \prox_{\tau \gamma_P f_P} \Big( \widetilde{\tdvar{Q}}^{k+1}\Big)\Bigg)$\\
\hspace*{0.4cm} 3. $\big(\widetilde{\tdvar{R}}^{k+1},\widetilde{\tdvar{S}}^{k+1}\big) = \sigma H(2 \widetilde{\tdvar{Q}}^{k+1} - \tdvar{Q}^{k}) + \big(\tdvar{R}^k, \tdvar{S}^k\big)$\\
\hspace*{0.4cm} 4. $\Big(\tdvar{R}^{k+1},\tdvar{S}^{k+1}\Big) = \Big(\widetilde{\tdvar{R}}^{k+1},\widetilde{\tdvar{S}}^{k+1}\Big) - \sigma\cdot \prox_{\gamma_{TV}/\sigma, \ell_1} \left(\frac{1}{\sigma} \widetilde{\tdvar{R}}^{k+1},\frac{1}{\sigma}\widetilde{\tdvar{S}}^{k+1}\right) $\\
\textbf{Stop if:}\\
\hspace*{0.4cm} $ \frac{\Vert \tdvar{Q}^{k+1} - \tdvar{Q}^{k}\Vert_2}{\Vert \tdvar{Q}^{k+1} \Vert_2} < 10^{-6}$ OR $k > 10^5$\\
\caption{Primal Dual algorithm}
\label{Algo:Prim-Dual}
\end{algorithm*}

This optimisation problem is solved by means of a primal-dual proximal algorithm, as in \cite{vu_splitting_2011, combettes_primal-dual_2011, condat_primal-dual_2013, komodakis_playing_2015}, which is particularly suited when the objective combines differentiable and non-differentiable functions along with linear operators. In such an iterative scheme, the non-differentiable functions are involved through their proximity operator \cite{moreau_proximite_1965} defined as:
\begin{equation}
(\forall u \in \mathcal{H}) \qquad \prox_{f}(u) = \operatorname{arg}\min_{x\in\mathcal{H}} f(x)+\frac{1}{2}\|u-x\|_2^2
\end{equation}
where $\mathcal{H}$ denotes a real Hilbert space and $f$ a convex, l.s.c., proper function from $\mathcal{H}$ to $]-\infty,+\infty]$. For further details on proximal algorithms, the reader could refer to \cite{combettes_proximal_2010, bauschke_convex_2011, parikh_proximal_2013}.

The proximity operator of the indicator of the convex set $C$ has a closed form expression as a projection \cite{theodoridis_adaptive_2011}: 
\begin{equation}
\prox_{\gamma_C f_C}(\tdvar{Q}) =\left\lbrace 
	\begin{array}{ll} 
	P_C(\tdvar{Q}) = \max(\tdvar{Q},\tdvar{B}) &\mbox{if } \gamma_C > 0\\
	\tdvar{Q} &\mbox{if } \gamma_C = 0.\\
	\end{array}
\right.
\end{equation}
The proximity operator of the function, $f_P$, also have a closed form expression \cite{combettes_douglas-rachford_2007}:
\begin{equation}
\begin{aligned}
&\prox_{\gamma_P f_P}(\tdvar{Q}) \\
&= \Big(\prox_{\gamma_P\psi}(B_{ij}^l,\eta^l Q_{ij}^l)\Big)_{(ijl)\in V\times V \times L} \\
&=  \Bigg(\frac{Q_{ij}^l - \gamma_P \eta^l  + \sqrt{\vert  Q_{ij}^l - \gamma_P \eta^l  \vert^2 + 4 \gamma_P B_{ij}^l}}{2}\Bigg)_{(ijl)\in V\times V \times L}.
\end{aligned}
\end{equation}

The proximal operator of the sum of these two functions satisfies the following property \cite{chaux_nested_2009}:
\begin{equation}
\prox_{\gamma_C f_C + \gamma_P f_P}(\tdvar{Q}) = P_C(\prox_{\gamma_P f_P}(\tdvar{Q})).
\end{equation}

Last, the $\ell_1$-norm, applied to $H$, as in Eq. \eqref{eqn:Obj:FTV}, also has a closed form expression for its proximity operator \cite{chaux_variational_2007, pustelnik_parallel_2009, donoho_-noising_1995, combettes_signal_2005}:
\begin{equation}
\begin{aligned}
\prox_{\gamma_{TV}\Vert \cdot \Vert_1}(\tdvar{R},\tdvar{S}) = \Bigg( & \sign(\tdvar{R})\ \max\lbrace \vert \tdvar{R}\vert - \gamma_{TV}, 0\rbrace,\\
 & \sign(\tdvar{S})\ \max\lbrace \vert \tdvar{S}\vert - \gamma_{TV}, 0\rbrace\Bigg) .
\end{aligned}
\end{equation}

The primal-dual proximal iterations designed for minimizing Eq. \eqref{eqn:Obj} are described in Algorithm \ref{Algo:Prim-Dual}. Under some technical assumptions regarding the domain of definition and the following condition \cite[theorem (3.1)]{condat_primal-dual_2013}:
\begin{equation}
\frac{1}{\tau} - \sigma \chi_H \geq \frac{\beta}{2},
\end{equation}
where the $\beta = \gamma_{TC} \beta_{TC} + \gamma_K \beta_K $ denotes the Lipschitz constant of $\gamma_{TC} f_{TC} + \gamma_K f_K$ and $\sigma>0$, the sequence $\big(\tdvar{Q}^{k+1}\big)_{k\in \mathbb{N}}$ converges to a minimizer of Eq. \eqref{eqn:Obj}.

\section{Simulated Case Study}
\label{sec:SCS}

	\subsection{Experimental setup}
		\subsubsection{Simulation context}
To test and validate the proposed method, a simplified road network model has been created. This has been preferred to a real case study for three reasons: tractability, the possibility to access the ground truth and the opportunity to explore the behaviour of the method for varied conditions. However, the number of nodes, the connectivity, the number of users and their OD patterns have been chosen to be consistent with those of a real networks.

The number of nodes of the simulated network is $\vert V\vert=50$ nodes. This number is kept relatively low to allow for a thorough exploration of the possible weights $\gamma_\cdot$ of problem \eqref{eqn:Obj}. For comparison, the Brisbane Bluetooth scanner network has around $900$ intersections equipped with vehicle identification devices. Other works on ODM estimation consider often few tens of nodes ($\propto 100$ OD flows) \cite{djukic_advanced_2015} while very recent works considered up to $300$ nodes \cite{perrakis_bayesian_2015}.

For the simulation, nodes are first located randomly on a grid and then links are created while aiming for an average connectivity of 6, a value consistent with that of real road networks \cite{porta_network_2006}. This is done first, by means of a minimum spanning tree (computed by the Kruskal's algorithm \cite{cormen_algorithms_2009}), then, by adding links randomly to the nodes with lower degree (sum of in and out edges) provided that the added links do not cross or repeat an existing one.

The number of users is fixed to $N = 10^5$. This leads to an mean flow per link of 3000 users. In big cities, it would correspond to around one hour of traffic during peak hours. 
Each node $i$ has a probability $p_O(i)$ of being an origin and similarly we define $p_D(j)$ the probability of node $j$ to be a destination.
Thus $\udvar{p_O}$ and $\udvar{p_D}$ satisfy:
\begin{equation}
	\label{eqn:ODProb}
	\left\lbrace
		\begin{array}{l}
		\udvar{D} = \esp(N\times \udvar{p_D}), \\
		\udvar{O} = \esp(N\times \udvar{p_O}). \\
		\end{array}
	\right.
\end{equation}
An origin and a destination are randomly associated to each user, according to the probabilities $p_O$ and $p_D$. We simulate a preferred direction of travel, consistent with the trends observed in urban context (mostly due to commuters). To this end, $p_O$ is decreasing linearly with the X-axis of the grid while $p_D$ is increasing linearly. The shortest path from origin to destination is then attributed to each user.

\begin{figure}[!h]
\begin{center}
	\includegraphics[width=8cm]{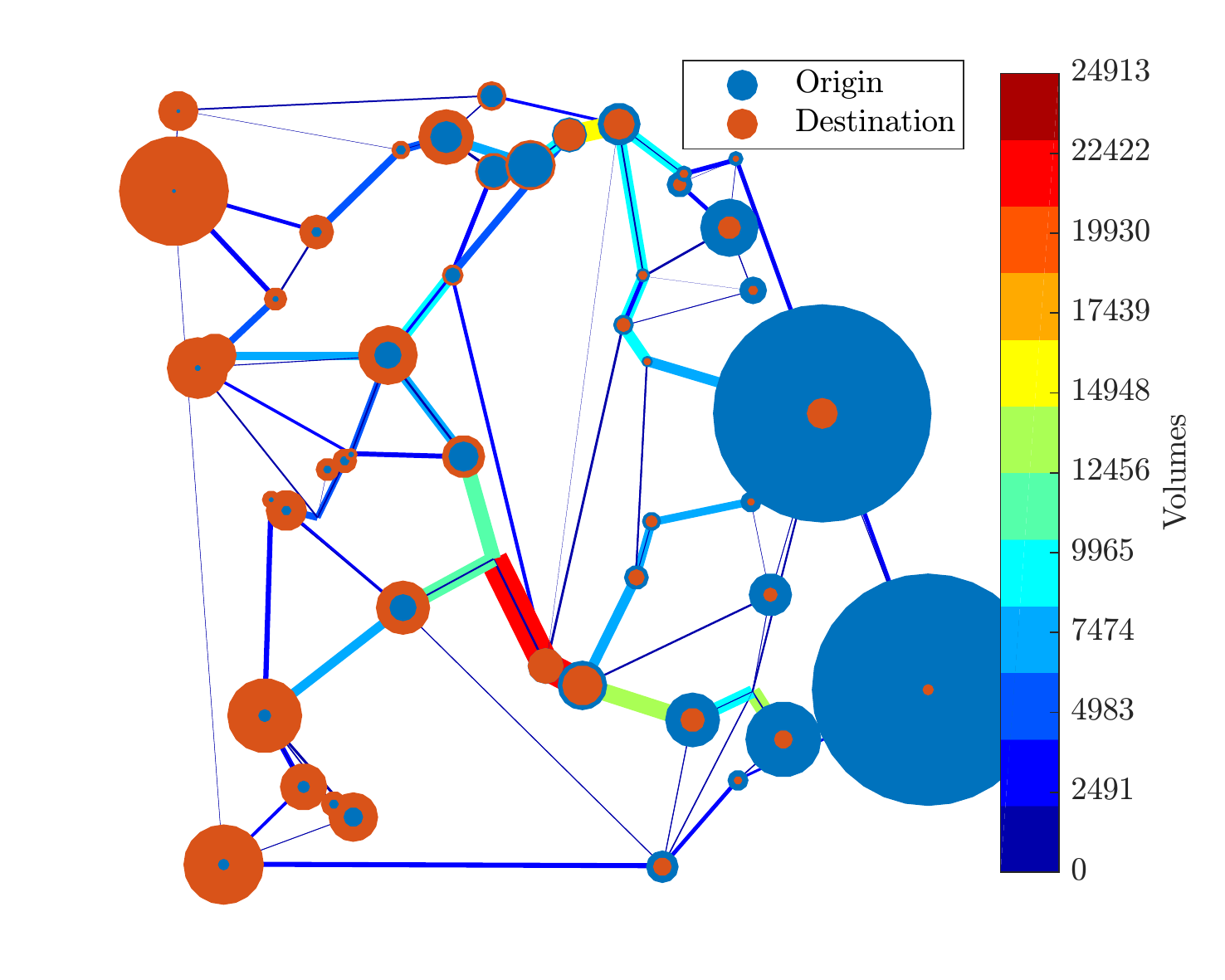}
	\caption{Simulated road networks with the projection of $\udvar{q}$ on the links: the width is proportional to the flows, also correlated with the color. For nodes, the color distinguishes between origin (blue) and destination flows (red) and the diameter of the nodes is proportional to their value in $\udvar{p_O}$ and $\udvar{p_D}$.}
	\label{fig:SimulatedNetwork}
\end{center}
\end{figure}

For each OD pair, a Bluetooth penetration rate is drawn from a Gaussian distribution of mean 30\% and standard deviation of 10\% (and truncated to be between 0 and 1). This choice accounts for the variability of the ownership distribution of Bluetooth devices (which is not known) from one node to another, depending, as an example, on the wealth of the neighbourhoods of the node. Each user has a probability equal to the Bluetooth penetration rate drawn for its OD of being equipped with a Bluetooth device. 
This gives us $\tdvar{B}$ while the full set of trajectories gives $\tdvar{Q}^\star$ for ground truth. The measured traffic flow per link $\udvar{q}$ is obtained from $\tdvar{Q}^\star$, assuming the addition a noise $\udvar{\varepsilon}$, for which each independent component is drawn from a Gaussian distribution ${\cal N}(0, r\cdot\udvar{q}^\star)$. For consistency with the noise usually measured on magnetic loops \cite{coifman_using_1999}, we take $r=5\%$.
Figure \ref{fig:SimulatedNetwork} illustrates the simulated case study with total volumes on the links ($\udvar{q}$) and the realisation of $p_O$ and $p_D$ for the $10^5$ users on the nodes.

		\subsubsection{Algorithmic parameters setup}
As discussed in section \ref{ssec:Algo}, the objective function \eqref{eqn:Obj} depends on five parameters $\gamma_{\cdot} \geq 0$. It appears however that, as $f_C$ can only be $0$ or $\infty$, exploring $\gamma_C \in \{0;1\}$ is enough. 
Moreover, the minimum in Equation \eqref{eqn:Obj} is preserved by a linear operation over the four remaining parameters $\gamma_{\cdot}$, and thus we choose $\gamma_P \in \{0;1\}$. This would then correspond to two situations: with or without considering the data stemming from sampled trajectories in the estimation process. We then explore the space of positive real numbers for the three remaining parameters $\gamma_{\cdot}$.

For those parameters, it has been observed that, for comparison purposes, it is justified to compare scenarii for rescaled values $\gamma_{\cdot} \beta_{\cdot}$. Indeed, $\beta_{\cdot}$ depends on the setup, in particular on $\udvar{q}$ and $\tdvar{B}$.

The algorithm stops if the convergence criteria is satisfied (\cf Alg. \ref{Algo:Prim-Dual}) or after $10^5$ iterations.

		\subsubsection{Performance evaluation}
The efficiency of the estimation algorithm is assessed by comparing its results to the ground truth $\tdvar{Q}^\star$ with two indicators. First, we denote RMSE the $\ell_2$ norm of the error divided by the norm of the ground truth (RMSE standing for Root Mean Square Error):
\begin{equation}
RMSE(\widehat{\tdvar{Q}}) = \frac{\Vert \widehat{\tdvar{Q}} - \tdvar{Q}^\star\Vert}{\Vert\tdvar{Q}^\star\Vert}.
\end{equation}
Second, EMD refers to the Earth Movers' Distance \cite{rubner_earth_2000}, a metric often used for image or distribution comparisons. It corresponds to the minimal cost that must be paid to transform the histogram of one image or distribution into the other.

We also provide a comparison of the estimates with two naive solutions  denoted $\widehat{\tdvar{Q}}^0$ and $\widehat{\tdvar{Q}}^1$, computed as the Bluetooth LODM multiplied by the mean Bluetooth penetration rate over the whole network ($\widehat{\eta}$), or over each link ($\udvar{\eta}$) respectively, for every $(i,j,l) \in  V \times  V \times  L$,
\begin{equation}
\label{eqn:Q0Q1}
\begin{aligned}
 & ( \widehat{{Q}}^0)^l_{ij}  = \widehat{\eta}\  B^l_{ij} &\mbox{where}\quad& \widehat{\eta} = \sum_l{\udvar{q}}/\sum_{i,j,l}\tdvar{B}, \\
 & ( \widehat{{Q}}^1)^l_{ij}  = \eta^l\  B^l_{ij} &\mbox{where}\quad&  \eta^l = q^l/\sum_{i,j,\bullet}\tdvar{B}.
\end{aligned}
\end{equation}

	\subsection{Results}
	\label{ssec:results}
		\subsubsection{Finding the best estimates}
\begin{figure}[t]
\begin{center}
	\includegraphics[width=8cm]{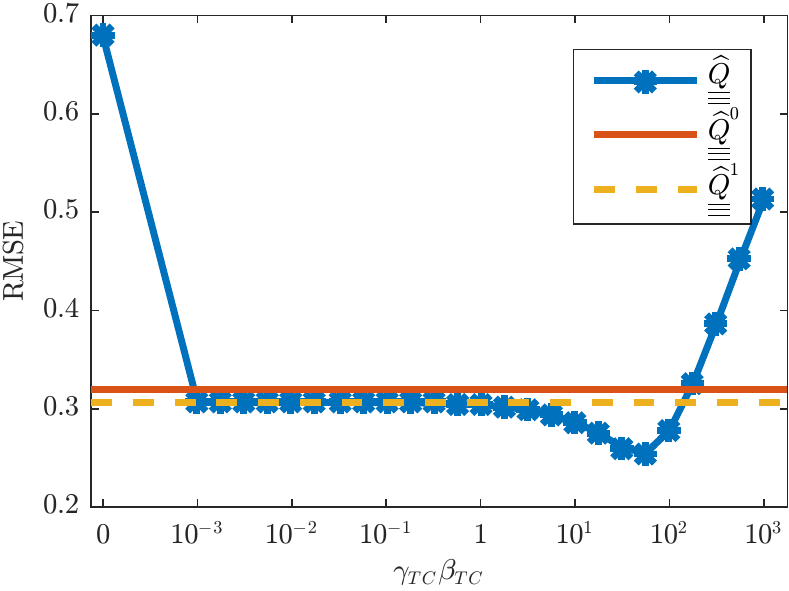}
	\caption{RMSE as a function of $\gamma_{TC}\ \beta_{TC}$ for $\gamma_K = \gamma_{TV} = 0$ and $\gamma_C = \gamma_P = 1$. Note that the first point on the left is for $\gamma_{TC}\ \beta_{TC} = 0$, to be distinguished from the logarithmic scale from the second point and after.}
	\label{fig:RMSE_Example}
\end{center}

\end{figure}

\begin{figure}[t]
\begin{center}
	\includegraphics[width=8cm]{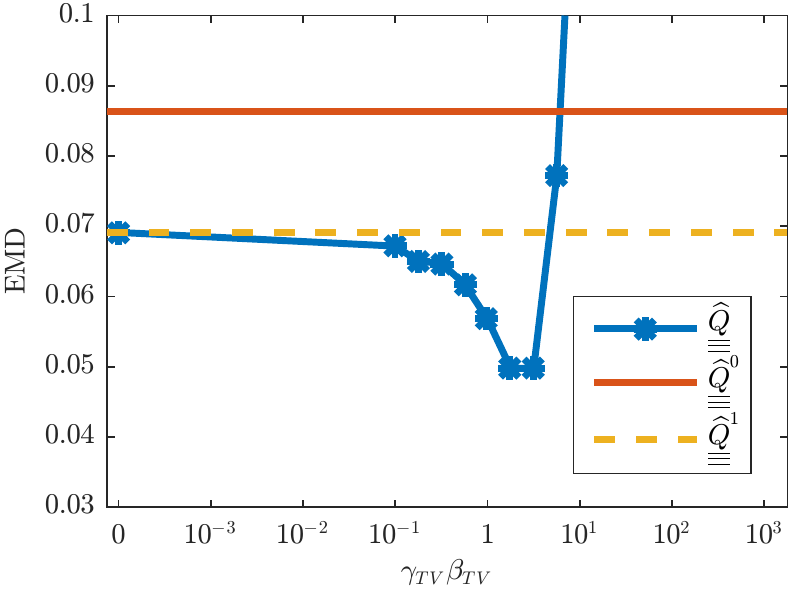}
	\caption{EMD as a function of $\gamma_{TV}\ \beta_{TV}$ for $\gamma_K = 0$ and $\gamma_{TC} = \gamma_C = \gamma_P = 1$. Again, note that the first point on the left is for $\gamma_{TV}\ \beta_{TV} = 0$ whereas the rest is drawn on a logarithmic scale.}
	\label{fig:EMD_Example}
\end{center}
\end{figure}

Solutions to problem \eqref{eqn:Obj} have been explored through a systematic exploration of the $\gamma_\cdot$ values within the positive real numbers.
Figures \ref{fig:RMSE_Example} and \ref{fig:EMD_Example} illustrate the evolution of the criteria RMSE (Fig. \ref{fig:RMSE_Example}) and EMD (Fig. \ref{fig:EMD_Example}) as a function of one $\gamma_\cdot$, the others being fixed. It highlights the existence of sets of parameters $\gamma_\cdot$ for which the estimates have minimal criteria while being lower than the criteria of the naive estimates $\widehat{\tdvar{Q}}^0$ and $\widehat{\tdvar{Q}}^1$. Therefore, in the following we denote by $\widehat{\tdvar{Q}}_{RMSE}$ the estimate $\widehat{\tdvar{Q}}$ minimizing the RMSE and by $\widehat{\tdvar{Q}}_{EMD}$ the one minimising the EMD. That is:
\begin{align}
  \widehat{\tdvar{Q}}_{RMSE}   \in & \underset{\widehat{\tdvar{Q}}}{\Argmin}  RMSE(\widehat{\tdvar{Q}}),\\
  \widehat{\tdvar{Q}}_{EMD} \in & \underset{\widehat{\tdvar{Q}}}{\Argmin} EMD(\widehat{\tdvar{Q}}).
\end{align}

\begin{figure*}[t]
	\includegraphics[width=17cm]{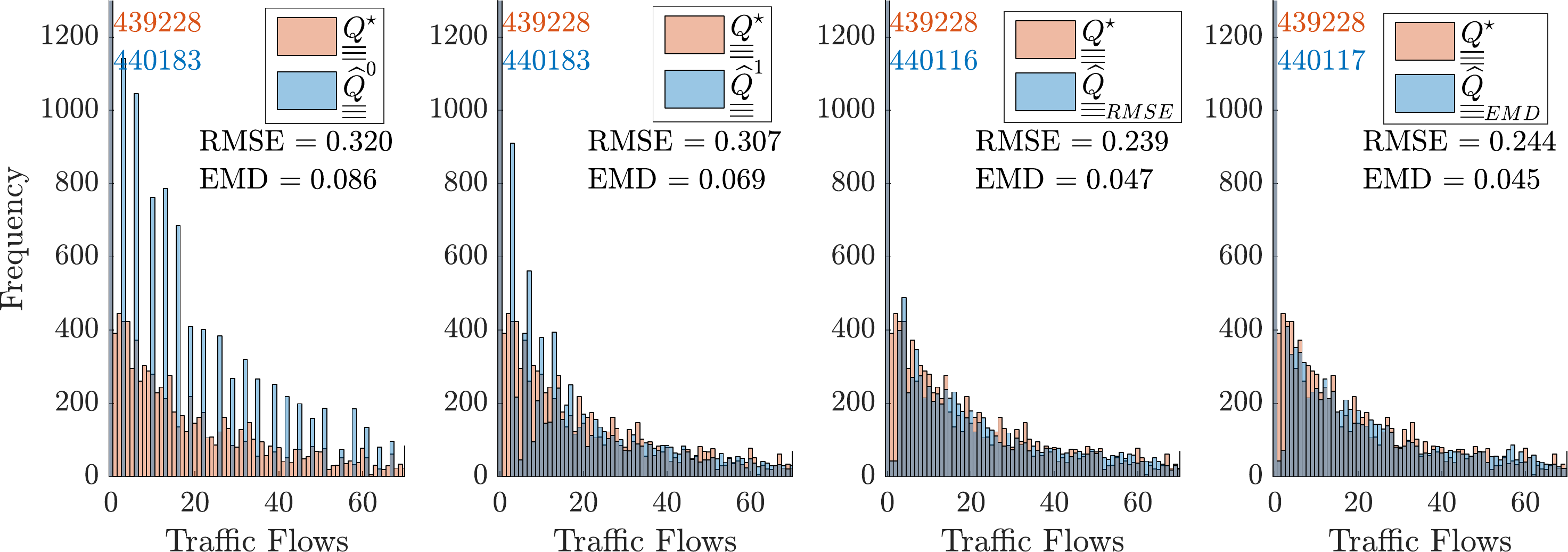}
	\caption{Distribution as histograms of the elements in $\widehat{\tdvar{Q}}^0$, $\widehat{\tdvar{Q}}^1$, $\widehat{\tdvar{Q}}_{RMSE}$, $\widehat{\tdvar{Q}}_{EMD}$ (in blue) superposed to the ground truth $\tdvar{Q}^\star$ (in red).}
	\label{fig:Histo}
\end{figure*}

Figure \ref{fig:Histo} represents the distribution of the elements in $\widehat{\tdvar{Q}}^0$, $\widehat{\tdvar{Q}}^1$, $\widehat{\tdvar{Q}}_{RMSE}$, $\widehat{\tdvar{Q}}_{EMD}$ as histograms, superposed to the ground truth $\tdvar{Q}^\star$ for comparison. As one could have expected from Eq \eqref{eqn:Q0Q1}, $\widehat{\tdvar{Q}}^0$ contains only multiples of the penetration rate value $\widehat{\eta}$. It has therefore a higher EMD value than other estimates.

	For those four same estimates, Table \ref{tbl:CompSol} presents the values for RMSE, EMD, $f_{TC}$ (consistency with observed counts) and $f_K$ (conformity with Kirchhoff's law). These two functions are chosen because they are the most important from a transport perspective point-of-view. When applicable, the corresponding values for the $\gamma_{\cdot}$ are indicated. The table illustrates that $\widehat{\tdvar{Q}}^0$ is a good solution from a transport perspective as it has low values for $f_{TC}$ and $f_K$ but is the worst performing when it is compared to the ground truth $\tdvar{Q}^\star$ (both for EMD and RMSE). 
	On the opposite, $\widehat{\tdvar{Q}}^1$ has better RMSE and EMD values than $\widehat{\tdvar{Q}}^0$ but performs poorly on the relevant transport indicators. Hence it is justified to actually solve problem \eqref{eqn:Obj} as none of the naive solutions bring satisfying results. While doing so, both $\widehat{\tdvar{Q}}_{RMSE}$ and $\widehat{\tdvar{Q}}_{EMD}$ give similar results in term of comparison to ground truth. 
	As one might expect, their performances with respect to the transport indicators are consistent with the weights applied to the corresponding functions: The solution $\widehat{\tdvar{Q}}_{RMSE}$ is reached for higher $\gamma_{TC}$ and, therefore, performs much better on the $f_{TC}$ criterion while $\widehat{\tdvar{Q}}_{EMD}$ has a relatively higher $\gamma_{K}$ and hence satisfies the Kirchhoff's law better.
	The question of which solution is the best partly depends on the reliability of the link counts (if the noise is high, satisfying perfectly $f_{TC}$ might not be relevant). Anyway, choosing between the two amounts to choosing the best $\gamma_\cdot$, a question left for future work.


\begin{table}[h]
\centering
\caption{Results for LODM estimates: the two naive solution and the ones minimizing RMSE and the EMD}
\label{tbl:CompSol}
\begin{tabular}{|l|lll|llll|}
\hline
         & $\gamma_{TC}$ & $\gamma_{K}$ & $\gamma_{TV}$ & RMSE & EMD   & $f_K$ & $f_{TC}$ \\
\hline
\hline
$\widehat{\tdvar{Q}}^0$       &             &            &             & 0.320 			& 0.086 		& \textbf{0} 	& 55   \\
$\widehat{\tdvar{Q}}^1$       &             &            &             & 0.307 			& 0.069 		& 1142   		& 84   \\
$\widehat{\tdvar{Q}}_{RMSE}$    & 31.6        & 0.008      & 0.015       & \textbf{0.239} & 0.047 		& 289    		& \textbf{1}\\
$\widehat{\tdvar{Q}}_{EMD}$   & 1           & 0.025      & 0.027       & 0.244 			& \textbf{0.045}& 133   		& 28   \\
\hline
\end{tabular}
\end{table}

		\subsubsection{Impact of each objective}

	If it appears from these results that it is justified to solve problem \eqref{eqn:Obj}, the question of the importance of each function can be raised. To answer such a question, Tables \ref{tbl:OneTwoConst}(a) (resp. (b)) summarises the best RMSE values (resp. EMD) when only the objectives indexed by the rows and column are not set to zeros. 
	Thus diagonal elements correspond to a single term in the objective function while the four others are set to zero and non diagonal elements involve at most the two terms indexed by the row and column. For example, element (1,2) of Table \ref{tbl:OneTwoConst}(a) corresponds to the best value of RMSE achieved for $\gamma_{TV} = \gamma_{P} = \gamma_{C} = 0$ and values of $\gamma_{K}$ and $\gamma_{TC}$ evaluated on a grid.
	For those tables, light-grey cells correspond to estimates that could not outperform the naive estimates, and darker grey elements, cases for which Algorithm \ref{Algo:Prim-Dual} reached the $100 000$ steps limit without convergence.

	The observation of both tables (for RMSE and EMD) leads to the conclusion that neither term gives satisfactory results by itself. None of the diagonal elements outperform the naive estimates. 
	To improve on those values, one must involve the Poisson assumption and either the traffic counts function or the Kirchhoff's law. This means that one cannot obtain a good estimate of the traffic flows while there are not at least one term ensuring data fidelity along with a regularization term. The fact that the Poisson function seems to be the most important can be expected as $\tdvar{B}$ brings the most information and confirms the importance of probe trajectories to solve such traffic problem.
	
	Thus in a second step, the Poisson function has been imposed ($\gamma_P=1$) with either one or two extra functions and similar results are gathered in Tables \ref{tbl:ThreeFourConst}(a) and (b). In those tables, the indicator function corresponding to the projection on the convex set $f_C$ has also been imposed ($\gamma_C=1$) for three reasons: First its additional computational cost is negligible compared to the other functions, second, it accelerates the convergence speed of the algorithm by reducing the number of steps required while having little impact on the values of the criteria at convergence. Last, it corresponds to the weakest assumption of our model: that the total flow is greater than the measured probe trajectories.

	In this second scenario (with $\gamma_P = \gamma_C=1$ and one or two additional functions), any estimate performs better than the naive one but for the case where only the total variation (TV) is added. Depending on whether a minimum is sought for RMSE or for the EMD, it is either the pair Kirchhoff's law and TV, or the pair Traffic Counts and TV, that are the best suited to complement the Poisson assumption.
	In any case, best results, as summarised in Table \ref{tbl:CompSol}, are achieved when all functions are involved. However, these results might question the role of the TV: First one can not reach better results than with naive estimates if at least the Poisson assumption and another function (K or TC) are not also involved. Second, the additional computation cost in Algorithm \ref{Algo:Prim-Dual} caused by the realisation of $H$, its adjoint and the proximal of the $\ell_1$-norm multiplies by four the time needed for each iteration (convergence reached in $\sim$4 hours instead of $\sim$1 on a Core i7 laptop). However, the total variation bring a slight improvement to the results for both RMSE (0.262 to 0.239) and the EMD (0.067 to 0.045). This is probably caused by the difficulty of estimating flows not sampled at all by the probe trajectories, task to which, in this version of the problem, the total variation is the only answer.

\begin{table}[h]
\centering
\caption{Best values of RMSE and EMD when only one or two constraints are considered}
\label{tbl:OneTwoConst}
\begin{tabular}{cc}
(a) & \begin{tabular}{|l|lllll|}
\hline
RMSE & TC    & K     & TV    & P     & C     \\
\hline
\hline
TC    & \cellcolor{lgray} 0.99 & \cellcolor{lgray} 0.98 	& \cellcolor{lgray} 0.99	& 0.30 							& \cellcolor{lgray} 0.67 \\
K     & 	  					& \cellcolor{lgray} 1 		& \cellcolor{lgray} 1		& \textbf{0.29}					& \cellcolor{lgray} 0.68 \\
TV    & 	  					& 	  					  	& \cellcolor{lgray} 1		& \cellcolor{dgray} 0.39		& \cellcolor{lgray} 0.68 \\
P     & 	  					& 	  						& 	  						& \cellcolor{lgray} 1			& \cellcolor{lgray} 0.68 \\
C     & 	  					& 	  						& 	  						& 	  							& \cellcolor{lgray} 0.68\\
\hline
\end{tabular} \\
 & \\
(b)& 
\begin{tabular}{|l|lllll|}
\hline
EMD & TC    & K     & TV    & P     & C     \\
\hline
\hline
TC  & \cellcolor{lgray} 2.03	& \cellcolor{lgray} 1.20	& \cellcolor{lgray} 1.20	& \textbf{0.07}					& \cellcolor{lgray} 1.43 \\
K   & 							& \cellcolor{lgray} 1.20	& \cellcolor{lgray} 1.20	& \cellcolor{lgray} 0.08		& \cellcolor{lgray} 0.81 \\
TV  & 							& 							& \cellcolor{lgray} 1.20	& \cellcolor{dgray} 0.19		& \cellcolor{lgray} 0.78 \\
P   & 							& 							& 							& \cellcolor{lgray} 1.20		& \cellcolor{lgray} 0.81 \\
C   & 							& 							& 							&				 				& \cellcolor{lgray} 0.81 \\
\hline
\end{tabular}\\
\end{tabular}
\end{table}

\begin{table}[h]
\centering
\caption{Best values of RMSE and EMD when $\gamma_C = \gamma_P =1$ and with one or two additional constraints}
\label{tbl:ThreeFourConst}
(a) \hfil (b)\hfil \\
\vspace{0.1cm}
\begin{tabular}{|l|lll|}
\hline
RMSE & TC    & K     & TV    \\
\hline
\hline
TC    & 0.27  & 0.26 & 0.26 \\
K     & 	  & 0.27 & \textbf{0.25} \\
TV    & 	  & 	 & \cellcolor{lgray} 0.35 \\
\hline
\end{tabular}
\hfil
\begin{tabular}{|l|lll|}
\hline
EMD & TC    & K     & TV    \\
\hline
\hline
TC  & 0.067 			& 0.067 & \textbf{0.046}			\\
K   & 					& 0.069 & 0.069 					\\
TV  & 					& 		& \cellcolor{lgray}0.605	\\
\hline
\end{tabular}
\hfil
\end{table}

	\subsection{Results with fewer users on the networks}
Finally, one might wonder whether this method, as presented in this article, can achieve similar results over smaller time periods, that is in our case, a lower number of users. In this section the main results are presented again for $N=10^4$ users. This correspond to an average flow of 300 vehicles per link. In a big city, this correspond to 5 to 10 minutes of traffic during peak hours. With such a low number of users we are reaching the limits of our model as $10^4$ users corresponds to $\sim$4 users (that is, in average, 1.3 probe trajectories) per OD. Therefore the impact of the Poisson assumption, inferring information from $\tdvar{B}$ decreases.
Table \ref{tbl:CompSol:10000} summarises the results in this case.
Yet there is still a 14\% improvement on the RMSE and a 30\% improvement on the EMD. These results are very encouraging as even in the limit cases, the estimates achieved with the algorithm are still an improvement with respect to the naive estimates.
\begin{table}[h]
\centering
\caption{Best Achieved Results with N=10 000}
\label{tbl:CompSol:10000}
\begin{tabular}{|l|lll|llll|}
\hline
                            & $\gamma_{TC}$ & $\gamma_{K}$ & $\gamma_{TV}$ & RMSE & EMD   & $f_K$ & $f_{TC}$ \\
\hline
\hline
$\widehat{\tdvar{Q}}^0$     &               &              &               & 0.398 			& 0.021 		 & \textbf{0}   & 13.3  		\\
$\widehat{\tdvar{Q}}^1$     &               &              &               & 0.396 			& 0.017 		 & 128   	    & \textbf{0.1}  \\
$\widehat{\tdvar{Q}}_{RMSE}$   & 1.78          & 0.25         & 0.027         & \textbf{0.341} & 0.013 		 & 10.7     	& 9.8   		\\
$\widehat{\tdvar{Q}}_{EMD}$ & 1             & 0.45         & 0.026         & 0.342 			& \textbf{0.012} & 5.8    	    & 13.7  			\\
\hline
\end{tabular}
\end{table}

\section{Conclusion}

We have shown that the Link dependent Origin-Destination matrix is an interesting tool for traffic representation. Moreover we have evidenced that its estimation can be performed with a primal dual algorithm and that the objective function to be minimized can be partially derived from natural properties of the problem (the consistency between measured and estimated traffic counts, the domain of definition and the Kirchhoff's law). Then by adding a few sensible relationships, as for example, the Poisson sampling assumption and the similarities between nearby flows computed as the total variation, one can obtain from such method, estimates that outperform the naive solutions. 
However improvements can still be sought, especially by designing new functions. 
Future works will demonstrate that, if available, traffic turn fractions at intersections can be involved as an additional term to achieve even better results. 
Another trail for developing this problem is to look for an online algorithm as the one presented in \cite[Section 5.2]{combettes_stochastic_2015}. Last, one could think about implementing time dependencies. This could be done by adding new relationships that would link successive estimations of the LODM or, alternatively, by using other methods: for example, a Kalman filter similarly to what have been done on traffic counts based ODM estimation \cite{cremer_new_1987}, or also with supplementary data (\eg Bluetooth \cite{barcelo_robustness_2013} or other sensors \cite{lu_kalman_2015}).

\printbibliography

\end{document}